\documentclass[12pt]{amsart}

\input xy
\xyoption{all}
\usepackage[english]{babel}
\usepackage[latin1]{inputenc}
\usepackage{graphicx}
\usepackage{amsmath}
\usepackage{amssymb}
\usepackage{amscd}
\usepackage{color}
\usepackage{oldgerm}
\usepackage{amsfonts}
\usepackage{newlfont}
\usepackage{longtable}
\usepackage{multirow}

\DeclareOldFontCommand{\rm}{\normalfont\rmfamily}{\mathrm}

\usepackage[all]{xy}

\usepackage{upref}

\def\C{\Bbb C}

\def\F{\Bbb F}

\def\ad{\operatorname{ad}}

\def\Aut{\operatorname{Aut}}

\def\Der{\operatorname{Der}}
\def\det{\operatorname{det}}

\def\dim{\operatorname{dim}}

\def\End{\operatorname{End}}

\def\Tr{\operatorname{Tr}}

\def\Ker{\operatorname{Ker}}

\def\Id{\operatorname{Id}}

\def\Rad{\operatorname{Rad}}

\def\Im{\operatorname{Im}}

\def\g{\frak g}
\def\gl{\frak{gl}}
\def\h{\frak h}

\def\s{\frak{sl}}

\theoremstyle{plain}\swapnumbers

\newtheorem{theorem}{Theorem}[section]

\newtheorem{Prop}[theorem]{Proposition}
\newtheorem{Def}[theorem]{Definition}
\newtheorem{Cor}[theorem]{Corollary}

\newtheorem{Remark}[theorem]{Remark}

\title[Generalized derivations]{Generalized derivations and Hom-Lie algebra structures on $\mathfrak{sl}_2$}
\author{R. Garc\'ia-Delgado}
\address{CIMAT-Unidad M\'erida}
\email{rosendo.garcia@cimat.mx}
\begin{document}
\maketitle
\begin{abstract}
The purpose of this paper is to show that there are Hom-Lie algebra structures on $\mathfrak{sl}_2(\mathbb{F}) \oplus \mathbb{F}D$, where $D$ is a special type of generalized derivation of $\mathfrak{sl}_2(\mathbb{F})$, and $\mathbb{F}$ is an algebraically closed field of characteristic zero. We study the representation theory of Hom-Lie algebras
within the approriate category and prove that any finite dimensional representation of a Hom-Lie algebra
of the form $\mathfrak{sl}_2(\mathbb{F}) \oplus \mathbb{F}D$, is completely reducible, in analogy to the well known Theorem of Weyl from the classical Lie theory. 
We apply this result to characterize the non-solvable Lie algebras 
of this type having an invertible $D$. It is shown that the generalized derivations $D$ of
$\mathfrak{sl}_2(\mathbb{F})$ that we study in this work, satisfy the Hom-Lie Jacobi identity for the Lie bracket of $\mathfrak{sl}_2(\mathbb{F})$.
Moreover, using root space decomposition techniques we provide an
intrinsic proof of the fact  that $\mathfrak{sl}_2(\mathbb{F})$ is the only simple Lie algebra
admitting non-trivial Hom-Lie structures.

\end{abstract}
\section{Introduction}

Extensions of Lie algebras can be regarded as mechanisms
to enlarge Lie algebras under certain prescriptions.
Some of the most studied and best understood extensions
are the semidirect products obtained from a given
Lie algebra $\g$ and a given derivation $D:\g\to\g$ of it.
At the end one obtains a new Lie algebra $\g[D]:=\g\oplus \mathbb{F} D$,
where $\mathbb{F} D$ stands for the one-dimensional
subspace of $\g[D]$ generated by $D$.
Then $\g$ becomes a Lie subalgebra of $\g[D]$ ---an ideal, actually---
and  the Lie bracket $[D,x]$ in $\g[D]$, is 
equal to $D(x)$, for any $x\in\g$.

\smallskip
\noindent
A generalized notion of derivation was introduced in [1]
and some special types of it were considered in [3] and [4].
See {\bf \S 2} below for a review.
We address
the question of studying extensions $\g[D]$ of the Lie algebra 
$\g=\frak{sl}_2(\Bbb F)$, over an algebraically closed field $\mathbb{F}$ of characteristic zero,
by some of these special types of generalized derivations
of the Lie product $[\,\cdot\,,\,\cdot\,]_\g$ of $\g$.
We shall keep denoting such generalized derivations by $D:\g\to\g$.
It turns out that the resulting extensions $\g[D]$ can be equipped with
well defined skew-symmetric, bilinear products that satisfy
a kind of {\it generalized Jacobi identity\/.} The algebraic 
structures thus obtained have appeared in the literature as examples
of what are now called Hom-Lie algebras.
Hom-Lie algebras
were apparently introduced in \cite{Hartwig} by considering
deformations of Lie algebras, and they are now being studied
for general skew-symmetric, bilinear, non-associative products
(see, for example \cite{Alejandra}, \cite{GSV}, \cite{Xie}, etc).

\smallskip
\noindent
Our aim is to classify, up to isomorphism, within the category of Hom-Lie algebras,
all the extensions $\mathfrak{sl}_2[D]
=\mathfrak{sl}_2(\mathbb{F}) \oplus \mathbb{F}D$,
over an algebraically closed field $\F$ of characteristic zero, by the special 
generalized derivations $D$ specified in \S 2 below.
As a matter of fact, each of these $D$'s satisfy the
Hom-Lie Jacobi identity for the Lie product of $\frak{sl}_2(\Bbb F)$.
It is proved in {\bf \S 4} below that the vector subspace of $\operatorname{End}_{\Bbb F}\left(\frak{sl}_2(\Bbb F)\right)$,
consisting of those linear maps $D$ that satisfy,
$$
\text{(HLJ)}\quad\quad\ \  \left[\,D(x), [y,z]\,\right]+\left[\,D(z), [x,y]\,\right]+\left[\,D(y), [z,x]\,\right]=0,
$$
can be decomposed in terms of two irreducible $\frak{sl}_2(\Bbb F)$-modules,
one of which is just the scalar multiples of the identity map
that obviously satisfy the ordinary Jacobi identity. The other
irreducible submodule is $5$-dimensional and consists exactly
of the generalized derivations of the specified type. 
Then, by making explicit
use of the appropriate canonical forms for the linear maps $D$ 
satisfying (HLJ) obtained in {\bf Prop \ref{forma canonica}},
we determine the isomorphism classes of the 
Hom-Lie algebras $\mathfrak{sl}_2[D]$ in {\bf Thm. \ref{theorem clasificacion}}.

\smallskip
\noindent
Once this classification is settled down, we may address the question of
determining the irreducible Hom-Lie modules for the Hom-Lie algebras $\mathfrak{sl}_2[D]$.
It turns out that these can be studied and characterized by using the 
representation theory of $\frak{sl}_2(\Bbb F)$. 
We prove in {\bf Thm. \ref{theorem-rep}} that any finite dimensional representation of a Hom-Lie algebra
of the form $\mathfrak{sl}_2(\mathbb{F}) \oplus \mathbb{F}D$, is completely reducible,
thus extending to the Hom-Lie category
the well known Theorem of Weyl from classical Lie theory.

\smallskip
\noindent
There are important differences between Lie algebras
admitting invertible derivations in the usual sense and Lie algebras
having invertible {\it generalized derivations\/.} For instance, a classical result,
due to Jacobson (see \cite{Jacobson}) states that any finite dimensional Lie
algebra having an invertible derivation must be nilpotent. 
As for a possible converse, Jacobson also raised the question
of whether an arbitrary nilpotent Lie algebra had an invertible derivation
defined on it or not. The answer to this question turned out to be negative
and a characterization of those nilpotent Lie algebras that cannot have
an invertible derivation was given in \cite{Dixmier}.
By way of significant contrast, it follows from {\bf Thm. \ref{theorem clasificacion}}
below, that $\mathfrak{sl}_2$ always admits invertible generalized derivations
of some special type. Furthermore, using the result that describes the Hom-Lie modules
of the Hom-Lie algebras $\mathfrak{sl}_2[D]$ (see {\bf Thm. \ref{theorem-rep}}),
we deduce a structure theorem for non-solvable Lie algebras admitting
a specific type of invertible generalized derivations (see {\bf Thm. \ref{no solubles}}).

\smallskip
\noindent
It was interesting for us to realize that any generalized derivation of $\frak{sl}_2(\Bbb F)$
of the type described in \S 2 below, satisfies $\text{(HLJ)}$. It is then natural to try 
to understand if something similar occurs in more general cases; namely, if
there are linear maps $D:\g\to\g$, not necessarily generalized derivations, 
that satisfy $\text{(HLJ)}$ for a simple Lie algebra $\frak{g}$.
We found that this problem was already solved in the literature.
As a matter of fact, it is proved in \cite{Xie} that the only linear maps $\g\to\g$
that satisfy  $\text{(HLJ)}$ for a simple Lie algebra $\g$ defined over 
an algebraically closed field of characteristic zero are the scalar mutiples of
the identity map as soon as the rank of $\g$ gets strictly bigger than $1$
(see {\bf Thm. 3.3} in \cite{Xie}). 
The proof provided in  \cite{Xie}, however, is an exhaustive case-by-case
trial with the aid of the computer package GAP (see https://www.gap-system.org).
Searching for a deeper understanding
of the problem, we found that the Hom-Lie problem for a simple Lie algebra $\g$
was also addressed in \cite{Jin}. They consider, however,  linear maps $\g\to\g$ 
that satisfy $\text{(HLJ)}$ and preserve the Lie product $[\,\cdot\,,\,\cdot\,]_\g$
(see {\bf Prop 2.2} in \cite{Jin}).
Using only root space decomposition techniques, we prove in {\bf \S 4} below, 
that for any linear map $\g\to\g$ ---not necessarily one that preserves
the Lie product in $\g$--- the Hom-Lie Jacobi identity is only 
satisfied for scalar multiples of the identity map whenever the rank
of the simple Lie algebra is greater than $1$. 
Therefore, it is only $\frak{sl}_2(\F)$ the
simple Lie algebra that can be equipped with linear maps
$\frak{sl}_2(\F)\to\frak{sl}_2(\F)$, other than multiples of the identity map,
that turn $\frak{sl}_2(\F)$ into a Hom-Lie algebra.

\section{Basic definitions and generalized derivations}

\begin{Def}{\rm
A Lie algebra over a field $\mathbb{F}$, is a vector space $\mathfrak{g}$ over $\mathbb{F}$, equipped with a skew-symmetric bilinear map 
$[\,\cdot\,,\,\cdot\,]:\mathfrak{g} \times \mathfrak{g} \rightarrow \mathfrak{g}$, 
satisfying the {\bf Jacobi identity}:
$$
[x,[y,z]]+[y,[z,x]]+[z,[x,y]]=0,\ \ \forall \, x,y,z \in \g.
$$
The element $[x,y]\in\g$ is referred to as 
{\it the Lie bracket\/} (or {\it the Lie product\/}) of $x$ and $y$ in $\g$.}
\end{Def}

\begin{Def}{\rm
Let $\mathfrak{g}$ be a Lie algebra over $\mathbb{F}$ and 
let $a,b,c \in \mathbb{F}$. A linear map $D:\g\to\g$ is a
\textbf{generalized derivation of type $(a,b,c)$ of} $\g$
(or an \textbf{$(a,b,c)$-derivation of} $\g$ for short), if
$$
a\,D([x,y])=b\,[D(x),y]+c\,[x,D(y)],\quad \forall\, x,y \in \mathfrak{g}.
$$ 
We shall denote by $\Der_{(a,b,c)}(\mathfrak{g})$ 
the vector subspace of $\operatorname{End}_{\F}(\frak{g})$
consisting of 
all generalized derivations of type $(a,b,c)$.}
\end{Def}

\noindent 
{\bf Note:} A derivation $D$ of a given Lie algebra $\g$, in its 
usual sense, is just a $(1,1,1)$-derivation.

\noindent The following result reduces the study of generalized derivations
to a few special cases (for more details see {\bf Thm. 1.1} in \cite{Dorado}, 
{\bf Thm. 2.2.1} in \cite{Hrivnak}, or {\bf Thm. 2.2} in \cite{H-Nov}.)

\begin{theorem}[\cite{Dorado}, Theorem 1.1]{\sl
Let $a,b,c\in\C$. Then $\Der_{(a,b,c)}(\g)$ can be obtained by 
intersection 
of two of the following subspaces with appropriate choices of the parameter $d$:
\begin{enumerate}
\item $\Der_{(d,0,0)}(\g)$,
\item $\Der_{(0,1,-1)}(\g)$,
\item $\Der_{(d,1,1)}(\g)$.
\end{enumerate}
}
\end{theorem}

\noindent 
{\bf Convention:}
In this work we shall deal with a particular case of generalized derivations;
namely, those of type $(a,1,1)$, corresponding to the case (3) of the above theorem. 
\smallskip

\noindent 
For simple Lie algebras and $(a,1,1)$-derivations, one has the following result
(see \cite{Burde2}, {\bf Thm. 5.11}):

\begin{theorem}[\cite{Burde2}, Theorem 5.11]\label{Burde}{\sl
Let $\mathfrak{g}$ be a simple Lie algebra over an algebraic closed field $\mathbb{F}$. Let $a \in \mathbb{F}$ be different from $-1$, $0$, $1$ or $2$. Then, $\Der_{(a,1,1)}(\mathfrak{g})=\{0\}$. Furthermore, $\Der_{(2,1,1)}(\mathfrak{g})=\mathbb{F}\operatorname{Id}_{\mathfrak{g}}$ and $\Der_{(1,1,1)}(\mathfrak{g})=\ad(\mathfrak{g})$.}
\end{theorem}
\begin{Remark}{\rm
The cases $\Der_{(0,1,1)}(\mathfrak{g})$ and $\Der_{(-1,1,1)}(\mathfrak{g})$
for a simple Lie algebra $\frak{g}$ are described in \S{\bf 4.1} and {\bf Thm. 3}, below.}
\end{Remark}

\section{Extensions by a generalized derivation}

\noindent Let $(\mathfrak{g},[\,\cdot\,,\,\cdot\,]_{\g})$ be a Lie algebra over $\mathbb{F}$ and let $D \in \Der(\mathfrak{g})$ be a derivation. The vector space $\g[D]=\mathfrak{g} \oplus \mathbb{F}D$, admits a Lie algebra structure, where the Lie bracket $[\cdot,\cdot]_{\g[D]}:\g[D] \times \g[D] \rightarrow \g[D]$, is defined by:
\begin{equation}\label{suma semidirecta 1}
\![x+\lambda D,y+\mu D]_{\g[D]}=[x,y]_{\mathfrak{g}}-\mu D(x)+\lambda D(y),\,\, \forall x,y \in\mathfrak{g},\,\, \forall \lambda,\mu \in \mathbb{F}.
\end{equation}
The Lie algebra $(\g[D],[\,\cdot\,,\,\cdot\,]_{\g[D]})$ is the semi-direct product of $(\mathfrak{g},[\,\cdot\,,\,\cdot\,]_{\g})$ and $\F\,D$. It is natural to ask what kind of algebraic structure
does the vector space  $\mathfrak{g} \oplus \mathbb{F} D$ have, when $D$
is an $(a,1,1)$-derivation.

\smallskip
\noindent 
Let $D \in \End_{\mathbb{F}}\mathfrak{g}$ be a non-zero $(a,1,1)$-derivation
and let $T:\mathfrak{g} \oplus \mathbb{F}D \rightarrow \mathfrak{g} \oplus \mathbb{F}D$, 
be the linear map defined by:
\begin{equation}\label{hom-lie map}
T(x+\lambda D)=x+a \lambda D;\ \ \quad  x \in \mathfrak{g},\quad \lambda \in \mathbb{F}.
\end{equation}
If $a\ne 1$, then the map $T$ cannot be the identity on $\g[D]$. It takes a straightforward
computation to verify that,
$$
[T(x),[y,z]_{\g[D]}]_{\g[D]}+[T(y),[z,x]_{\g[D]}]_{\g[D]}
+[T(z),[x,y]_{\g[D]}]_{\g[D]}=0,
$$
for all $x,y,z \in \g[D]$. Thus, the vector space $\g[D]=\mathfrak{g} \oplus \mathbb{F}D$, with the bracket $[\,\cdot\,,\,\cdot\,]_{\g[D]}$, defined 
on it, is not a Lie algebra if $a \neq 1$, but a \emph{Hom-Lie algebra} in the sense of the following:

\begin{Def}\label{hom-Lie Def}{\rm
Let $\g$ be a vector space over $\mathbb{F}$. Let $T \in \End_{\mathbb{F}}(\g)$ be a linear map and let $[\,\cdot\,,\,\cdot\,]:\g \times \g \rightarrow \g$ be a skew-symmetric bilinear map satisfying:
$$
[T(x),[y,z]_{\g}]_{\g}+[T(y),[z,x]_{\g}]_{\g}+[T(z),[x,y]_{\g}]_{\g}=0,\quad \forall x,y,z \in \g.
$$
The triple $(\g,[\,\cdot\,,\,\cdot\,]_{\g},T)$ is a \textbf{Hom-Lie algebra over  $\mathbb{F}$}.
}
\end{Def}

To simplify, throughout this work, we simply say that $\g[D]$, with the bracket defined by \eqref{suma semidirecta 1} and the linear map $T \in \End_{\F}(\g[D])$ as \eqref{hom-lie map}, is the Hom-Lie algebra defined by $D$ and $T$.
\subsection{Extensions of $\mathfrak{sl}_2(\F)$ by a $(-1,1,1)$-derivation}

Observe that the above Hom-Lie algebra extension makes sense {\it for any\/} 
Lie algebra admitting a non-zero $(a,1,1)$-derivation. 
Let $\mathfrak{g}$ be a simple Lie algebra. As we have already observed,
if $a=1$, the obtained extension is
none other than the usual semi-direct product of $\frak{g}$ by an inner derivation $D$ 
of $\mathfrak{g}$. In order to obtain something new, consider $a \neq 1$.  
Since $\mathfrak{g}$ is simple, {\bf Thm. \ref{Burde}} restricts the possibilities 
of $a$ to be either $0$, $-1$ or $2$. 
For $a=2$ we have $\operatorname{Der}_{(2,1,1)}(\mathfrak{g})=
\mathbb{F}\operatorname{Id}_{\mathfrak{g}}$ (see {\bf Thm. 2}). 
It remains to understand the cases $a=0$ and $a=-1$. 
It is proved in  {\bf Lemma 6.1} of \cite{Leger}, that for a simple Lie algebra $\frak{g}$,
$\operatorname{Der}_{(0,1,1)}(\mathfrak{g})=\{0\}$. Thus, the case $a=0$ is trivial.
This leaves us with the case $a=-1$. 
Now, the proof of the following result can be found in \cite{Burde2} (see {\bf Thm. 5.12}):

\begin{theorem}[\cite{Burde2}, Theorem 5.12]{\sl
Let $\mathfrak{g}$ be a simple Lie algebra of rank at least two. Then $\Der_{(-1,1,1)}(\mathfrak{g})=0$.}
\end{theorem}

\noindent 
Therefore, the only possibility for a simple Lie algebra $\frak{g}$ to admit non-zero
$(-1,1,1)$-derivations is that $\frak{g}=\mathfrak{sl}_2(\F)$. 
Let $\{H,E,F\}$ be the standard basis for $\mathfrak{sl}_2(\F)$, so that,
$[H,E]=2E$, $[H,F]=-2F$ and $[E,F]=H$. 
It is shown in \cite{Burde2} that $\Der_{(-1,1,1)}(\mathfrak{sl}_2(\F))$
is the $5$-dimensional vector subspace of $\operatorname{End}(\frak{sl}_2)$
generated by the linear transformations whose 
$3 \times 3$ matrices in the basis $\{H,E,F\}$ are given by:
$$
P=\begin{pmatrix}
2 & 0 & 0\\
0 & -1 & 0\\
0 & 0 & -1
\end{pmatrix}\!,\quad
Q=\begin{pmatrix}
0 & 1 & 0\\
0 & 0 & 0\\
2 & 0 & 0
\end{pmatrix}\!,\quad
R=\begin{pmatrix}
0 & 0 & 1\\
2 & 0 & 0\\
0 & 0 & 0
\end{pmatrix}\!,
$$
\begin{equation}\label{base}
S=
\begin{pmatrix}
0 & 0 & 0\\
0 & 0 & 1\\
0 & 0 & 0
\end{pmatrix}\!,\quad
T=
\begin{pmatrix}
0 & 0 & 0\\
0 & 0 & 0\\
0 & 1 & 0
\end{pmatrix}\!.
\end{equation}
Thus, we
associate a $5$-tuple $(\zeta,\eta,\sigma,\lambda,\mu) \in \mathbb{F}^5$
to each $(-1,1,1)$-derivation $D$, where
$$
D=\zeta P+\eta Q+\sigma R+\lambda S+\mu T=\begin{pmatrix}
2\zeta & \eta & \sigma \\
2 \sigma & -\zeta & \lambda \\
2 \eta & \mu & -\zeta
\end{pmatrix}\!,
$$
From now on, we write $D=(\zeta,\eta,\sigma,\lambda,\mu)$ for such a $(-1,1,1)$-derivation $D$ of $\mathfrak{sl}_2(\F)$.

In order to classify the Hom-Lie algebras of the form
$\frak{sl}_2\oplus \mathbb{F} D$,
we take into account the following general definition (see \cite{Hartwig}, page 331):

\begin{Def}\label{Hom-Lie morfismo}{\rm
Let $(\mathfrak{g},[\cdot,\cdot]_{\mathfrak{g}},T)$ and 
$(\mathfrak{h},[\cdot,\cdot]_{\mathfrak{h}},S)$ 
be Hom-Lie algebras. A \textbf{Hom-Lie algebra morphism} 
from $(\mathfrak{g},[\cdot,\cdot]_{\mathfrak{g}},T)$ into 
$(\mathfrak{h},[\cdot,\cdot]_{\mathfrak{h}},S)$,
is a linear map, $\psi:\mathfrak{g} \rightarrow \mathfrak{h}$, satisfying,
\begin{itemize}

\item[(i)] $\psi([x,y]_{\mathfrak{g}})=[\psi(x),\psi(y)]_{\mathfrak{h}}$, for all $x,y \in \mathfrak{g}$,

\item[(ii)] $\psi \circ T=S \circ \psi$.

\end{itemize}
}
\end{Def}
In particular, for the family of Hom-Lie algebras $\mathfrak{sl}_2(\F)[D]$, obtained through
non-trivial $(a,1,1)$-derivations $D$, we have the following:

\begin{Prop}\label{criterio de isomorphismo}{\sl
Let $\Bbb F$ be a field with characteristic different from $2$.
Let $D,D^\prime \in \Der_{(-1,1,1)}(\mathfrak{sl}_2(\F))$ and let $T \in \End_{\mathbb{F}}(\mathfrak{sl}_2(\F)[D])$, defined by $T(x+\xi D)=x-\xi D$, for all $x \in \mathfrak{sl}_2$ and for all $\xi \in \F$, (resp $T^\prime \in \End_{\mathbb{F}}(\mathfrak{sl}_2(\F)[D]^\prime)$). Then, the Hom-Lie algebras $(\mathfrak{sl}_2(\F)[D],T)$ and $(\mathfrak{sl}_2(\F)[D]^\prime,T^\prime)$ are isomorphic if and only if there is an automorphism $g \in \operatorname{Aut}(\mathfrak{sl}_2(\F))$ and a non-zero scalar 
$\xi \in \mathbb{F}-\{0\}$, such that $D^{\prime}=\xi \, g\circ D \circ g^{-1}$.}
\end{Prop}
\begin{proof} 
Let$D$ and $D^\prime$ be two $(-1,1,1)$-derivations of $\mathfrak{sl}_2(\F)$ and let $\psi: \mathfrak{sl}_2(\F)[D] \rightarrow \mathfrak{sl}_2(\F)[D]^{\prime}$ be a Hom-Lie algebra isomorphism. Let $\pi_{\mathfrak{sl}_2(\F)}:\mathfrak{sl}_2(\F)[D] \to \mathfrak{sl}_2(\F)$ and $\pi_{D^{\prime}}:\mathfrak{sl}_2(\F)[D]\to \mathbb{F}D^{\prime}$ be the 
linear projections onto $\mathfrak{sl}_2(\F)$ and $\mathbb{F}D^{\prime}$, respectively. Let $g=
\pi_{\mathfrak{sl}_2(\F)} \circ \psi\vert_{\frak{sl}_2}:\mathfrak{sl}_2(\F) \to \mathfrak{sl}_2(\F)$ and let $\varphi:\mathfrak{sl}_2(\F) \to \mathbb{F}$ be the linear map for which $(\pi_{D^{\prime}} \circ \psi)(x)=\varphi(x)D^{\prime}$, for all $x \in \mathfrak{sl}_2(\F)$. 
Then, $\psi(x)=g(x)+\varphi(x)D^{\prime}$, for all $x \in \mathfrak{sl}_2(\F)$. Also, there exists $x' \in \mathfrak{sl}_2(\F)$ and $\xi' \in \mathbb{F}$, such that $\psi(D)=x'+\xi' D^{\prime}$. 
Take $x$ and $y$ in $\mathfrak{sl}_2(\F)$. Definition \ref{Hom-Lie morfismo}.(i) says that, 
$\psi([x,y])=g([x,y])+\varphi([x,y])D$. On the other hand
$$
\aligned
\,[\psi(x),\psi(y)]&=[g(x)+\varphi(x)D^{\prime},g(y)+\varphi(y)D^{\prime}]\\
\,& =[g(x),g(y)]-\varphi(y)D^{\prime}(g(x))+\varphi(x)D^{\prime}(g(y)).
\endaligned
$$
Since $\mathfrak{sl}_2(\F)=[\mathfrak{sl}_2(\F),\mathfrak{sl}_2(\F)]$, it follows that
$\varphi=0$ and $g \in \operatorname{Aut}(\mathfrak{sl}_2(\F))$.
On the other hand, Definition \ref{Hom-Lie morfismo}.(ii), states that $\psi \circ T(D)=T^\prime \circ \psi(D)$.
Therefore, $\psi(T(D))=\psi(-D)=-x'-\xi' D^{\prime}$, and $T^\prime(\psi(D))=T^\prime(x'+\xi'D^{\prime})=x'-\xi' D^{\prime}$, thus $x'=0$. This shows that $\psi(D)=\xi' D^{\prime}$. 
\smallskip

\noindent 
Now, let $x \in \mathfrak{sl}_2(\F)$. Since $\psi([D,x])=[\psi(D),\psi(x)]$, we deduce that
$g(D(x))=\xi' D^{\prime}(g(x))$. Therefore, $g\,\circ\,D=\xi^\prime\,D^{\prime}\,\circ\,g$, and 
$D^{\prime}=(\xi^\prime)^{-1}\, g\,\circ\,D\,\circ\, g^{-1}$. The converse statement is obvious.
\end{proof}

\begin{Remark}{\rm 
The proof of this Prop makes us realize that
there is a group action behind the isomorphism between two
extensions of $\frak{sl}_2$ by $(-1,1,1)$-derivations.
Namely, let $G=\mathbb{F}- \{0\} \times \operatorname{Aut}(\mathfrak{sl}_2(\F))$,
and consider the left $G$-action on $\operatorname{Der}_{(-1,1,1)}(\mathfrak{sl}_2(\F))$
given by,
$$
(\xi,g).D \mapsto \xi \, g \circ D \,\circ \,g^{-1},\quad \xi \in \mathbb{F}- \{0\}\ \text{and}\ g \in\operatorname{Aut}(\mathfrak{sl}_2(\F)).
$$
Therefore, $\mathfrak{sl}_2(\F)[D]$ and $\mathfrak{sl}_2(\F)[D]^\prime$ are isomorphic if and only if there exists a $(\xi,g) \in \mathbb{F}-\{0\} \times \operatorname{Aut}(\mathfrak{sl}_2(\F))$ such that $D^{\prime}=(\xi,g).D$. We shall write $D\sim D^\prime$ and we say that
$D$ and $D^\prime$ are equivalent if they belong to the same $G$-orbit.
Note, however, that any $D \in \operatorname{Der}_{(-1,1,1)}(\mathfrak{sl}_2(\F))$ is 
trivially equivalent to 
$D^\prime = \xi\,D$, for any $\xi \in \F - \{0\}$. 
It follows that $D$ and $D^\prime$ belong to the same $G$-orbit if and only 
if there exists an automorphism $g \in \Aut(\mathfrak{sl}_2(\F))$, 
such that $D^{\prime}=g\circ D \circ g^{-1}$.}
\end{Remark}

\subsubsection{The automorphism group of $\mathfrak{sl}_2(\F)$}

We shall now determine a set of canonical forms for the left $\operatorname{Aut}(\mathfrak{sl}_2)$-action
$g.D\mapsto g\cdot T= g^{-1} \circ D \circ g$ on $\operatorname{Der}_{(-1,1,1)}(\mathfrak{sl}_2)$.
It is known from \cite{Hum} that $\operatorname{Aut}(\mathfrak{sl}_2)$ is a semidirect product between $\operatorname{Inn}(\mathfrak{sl}_2)$ and the group of \emph{diagonal automorphism} $\Gamma(\mathfrak{sl}_2)$, \emph{i.e.}, those automorphism which are the identity on $H$ and a scalar multiplication in each root space. The group $\operatorname{Inn}(\mathfrak{sl}_2)$ is the subgroup of inner automorphisms
of $\mathfrak{sl}_2$; this in turn is determined by those elements
$x \in \mathfrak{sl}_2$ for which $\ad\vert_{\mathfrak{sl}_2}(x)$ is nilpotent;
these are known as $\ad\!\vert_{\mathfrak{sl}_2}$-{\it nilpotent elements\/,}
or $\ad$-{\it nilpotent\/} for short within our context.
It is known from \cite{Hum} (see
exercise $7$, chapter IV, page $88$) that $\operatorname{Aut}(\mathfrak{sl}_2(\F))
=\operatorname{Inn}(\mathfrak{sl}_2(\F))$, 
where $\operatorname{Inn}(\mathfrak{sl}_2(\F))$ is the subgroup of inner automorphisms
of $\mathfrak{sl}_2(\F)$; this in turn is determined by those elements
$x \in \mathfrak{sl}_2(\F)$ for which $\ad\vert_{\mathfrak{sl}_2(\F)}(x)$ is nilpotent;
these are known as $\ad\!\vert_{\mathfrak{sl}_2(\F)}$-{\it nilpotent elements\/,}
or $\ad$-{\it nilpotent\/} for short within our context.
It is easy to see that $x \in \mathfrak{sl}_2(\F)$ is $\ad$-nilpotent if and only if,
$$
x \in \mathbb{F}E,\quad x \in \mathbb{F}F,\quad x \in \{aH+acE-ac^{-1}F\,|\,a,c \in \mathbb{F}- \{0\}\}.
$$
Thus, $\operatorname{Aut}(\mathfrak{sl}_2(\F))$ is generated 
by the following automorphisms of $\frak{sl}_2$, 
$$
g_a=\exp(\ad\!\vert_{\mathfrak{sl}_2(\F)}(aE)),\quad  h_a=\exp(\ad\!\vert_{\mathfrak{sl}_2(\F)}(aF)),\quad
\text{and}
$$
$$
f_{a,c}=\exp(\ad\!\vert_{\mathfrak{sl}_2(\F)}(aH+acE-ac^{-1}F));\quad a,c \in \mathbb{F}- \{0\}.
$$
The matrices of $g_a$, $h_a$ and $f_{a,c}$ in the standard basis $\{H,E,F\}$, are:
$$
g_a=
\begin{pmatrix}
1 & 0 & a\\
-2a & 1 & -a^2\\
0 & 0 & 1
\end{pmatrix},\quad
h_a=
\begin{pmatrix}
1 & -a & 0\\
0 & 1 & 0\\
2a & -a^2 & 1
\end{pmatrix},\quad\text{and},
$$
$$
f_{a,c}=
\begin{pmatrix}
1-2a^2 & ac^{-1}+a^2c^{-1} & ac-a^2c\\
-2ac-2a^2c & (1+a)^2 & -a^2c^2\\
-2ac^{-1}+2a^2c^{-1} & -a^2c^{-2} & (1-a)^2
\end{pmatrix}\!,\quad a,c \in \mathbb{F}- \{0\}.
$$
Thus, the $\Aut(\mathfrak{sl}_2(\F))$-action $D \mapsto g^{-1} \circ D \circ g$ 
on $\operatorname{Der}_{(-1,1,1)}(\mathfrak{sl}_2(\F))$ can be written for the generators of
$\Aut(\mathfrak{sl}_2(\F))$ as,
$$
K_a(D)=g_a^{-1} \circ D \circ \, g_a,\quad L_a(D)=h_a^{-1} \circ D \circ \,h_a,
\quad J_{a,c}(D)=f_{a,c}^{-1}\circ D\circ \,f_{a,c}.
$$
Observe that $g_a^{-1}=g_{-a}$, $h_a^{-1}=h_{-a}$ and $f_{a,c}^{-1}=f_{-a,c}$, for all $a,c \in \mathbb{F}-\{0\}$. 
Now, using the identification of $\operatorname{Der}_{(-1,1,1)}(\mathfrak{sl}_2(\F))$ with $\Bbb F^5$, the images of a given 
$D=(\zeta,\eta,\sigma,\lambda,\mu)$, 
under $K_a$, $L_a$ and $J_{a,c}$, are respectively given by,
$$
\aligned
K_a(D)&=
(\zeta+2a \eta+a^2 \mu, \eta+a \mu,3a \zeta-3a^2 \eta+\sigma-a^3 \mu,
\\
&\qquad 
6a^2 \zeta+4a^3 \eta-4a \sigma+\lambda+a^4 \mu, \mu),\\
L_a(D)&=
(\zeta-2a \sigma+a^2 \lambda,3a \zeta+\eta-3a^2 \sigma+a^3 \lambda,\sigma-a \lambda,\lambda,
\\
&\qquad
6a^2 \zeta+4a \eta-4a^3 \sigma+a^4 \lambda+\mu),\\
J_{a,c}(D)&=
(J_{a,c}^1(D),J_{a,c}^2(D),J_{a,c}^3(D),J_{a,c}^4(D),J_{a,c}^5(D)),
\endaligned
$$
where the components $J_{a,c}^k:\mathbb{C}^5 \rightarrow \mathbb{C}$, $\,1 \leq k \leq 5\,$, 
are given by,
$$
\aligned
J_{a,c}^1(D)&=
(-1+6a^2-6a^4)\zeta-2ac(1-a)(1-2a^2)\eta
\\
&\quad 
-2ac^{-1}(1+a)(1-2a^2)\sigma
+a^2c^2(1+a)^2\lambda+a^2c^2(1-a)^2\mu
\\
J_{a,c}^2(D)&=
3ac^{-1}(1-a)(1+2a)\zeta+(1-a)^2(1-4a^2)\eta
\\
&\quad 
+a^2c^{-2}(4a^2-3)\sigma-a^3c^{-3}(1+a)\lambda+ac(1-a)^3\mu,
\\
J_{a,c}^3(D)&=
-3ac(1+a)(1+2a)\zeta +a^2c^2(4a^2-3)\eta
\\
&\quad
+(1+a)^2(1-4a^2)\sigma-ac^{-1}(1+a)^3\lambda-a^3c^3(1-a)\mu,
\\
J_{a,c}^4(D)&=
6a^2c^2(1+a)^2\zeta+4a^3c^3(1+a)\eta-4ac(1+a)^3\sigma
\\
&\quad
+(1+a)^4\lambda+a^4c^4\mu,\\
J_{a,c}^5(D)&=
6a^2c^{-2}(1-a)^2\zeta-4ac^{-1}(1-a)^3\eta+4a^3c^{-3}(1-a)\sigma
\\
&\quad
+a^4c^{-4}\lambda+(1-a)^4\mu.
\endaligned
$$

\subsubsection{Canonical forms for $(-1,1,1)$-derivations of $\mathfrak{sl}_2(\F)$}

In order to obtain a classification for those Hom-Lie algebras of the form $\mathfrak{sl}_2(\F)[D]$, we shall first obtain the canonical forms for the $(-1,1,1)$-derivations of $\mathfrak{sl}_2(\F)$. The group action is  $g.D=g \circ D \circ g^{-1}$, where $D \in \operatorname{Der}_{(-1,1,1)}(\mathfrak{sl}_2(\F))$ and $g \in \Aut(\mathfrak{sl}_2(\F))$.

\begin{Prop}\label{forma canonica}{\sl
Let $\mathbb{F}$ be an algebraically closed field of characteristic zero. Under the left
$\Aut(\mathfrak{sl}_2(\F))$-action on $\operatorname{Der}_{(-1,1,1)}(\mathfrak{sl}_2(\F))$,
$D\mapsto g.D=g\circ D \circ g^{-1}$, each $D$ admits a choice of $g$ such that,
$$
g.D=\begin{pmatrix}
0 & \nu^{-1}\eta & \nu\sigma \\
2 \nu \sigma & 0 & \nu^{2} \lambda \\
2 \nu^{-1}\eta & 0 & 0
\end{pmatrix}\!,\quad \eta,\sigma,\lambda \in \mathbb{F},\,\,\nu \in \F - \{0\}.
$$
}
\end{Prop}
\begin{proof}
\noindent Let $D\in \operatorname{Der}_{(-1,1,1)}(\mathfrak{sl}_2)$
and suppose $\mu \neq 0$. If $(\zeta,\eta,\sigma,\lambda) \neq (0,0,0,0)$, 
then $p=6\zeta x^2+4\eta x-4\sigma x^3+\lambda x^4+\mu$ is a non-constant 
polynomial with $p(0)=\mu\neq 0$. 
Thus, there exists $a \in \mathbb{F} -\{0\}$ such that $p(a)=0$. Applying $L_a$ to $D$, we get that $D$ is equivalent to $D^{\prime}=(\zeta^{\prime},\eta^{\prime},\sigma^{\prime},\lambda^{\prime},0)$.
On the other hand, if
$\zeta=\eta=\sigma=\lambda=0$, then
$D=(0,0,0,0,\mu)$,
and $J_{1,c}(0,0,0,0,\mu)=(0,0,0,c^{-4}\mu,0)$. 
In any case, 
$D$ is equivalent to 
$(\zeta^{\prime},\eta^{\prime},\sigma^{\prime},\lambda^{\prime},0)$. 
Therefore, with no loss of generality, we may always assume $\mu=0$
and start with $D$ of the form $D=(\zeta,\eta,\sigma,\lambda,0)$.
\smallskip

We now want to prove that any $D=(\zeta,\eta,\sigma,\lambda,0)$ is equivalent to $D^{\prime}=(0,\eta^{\prime},\sigma^{\prime},\lambda^{\prime},0)$. If $\zeta=0$, there is nothing to prove.
Thus, we assume from the start $\zeta \neq 0$.
If $\eta \neq 0$, we apply $K_{-\zeta/(2\eta)}$ to $D$, this
make $D$ equivalent to $D^\prime=(0,\eta^{\prime},\sigma^{\prime},\lambda^{\prime},0)$. 
So, let us assume $\eta=0$, \emph{i.e.,} $D=(\zeta,0,\sigma,\lambda,0)$, with $\zeta \neq 0$.
Now, suppose first $\sigma \neq 0$ or $\lambda \neq 0$
and take any $a^{\prime}\in\F - \{0\}$.  
Applying 
$K_{a^{\prime}}$ to $D$, makes the latter equivalent 
to $D^\prime=(\zeta^{\prime},0,\sigma^{\prime},\lambda^{\prime},0)$, where $\zeta^{\prime}=\zeta$ and:
$$
\sigma^{\prime}=-3a^{\prime}\zeta+\sigma,\quad \lambda^{\prime}=-6{a^{\prime}}^2\zeta-4a^{\prime}\sigma+\lambda,\quad a^{\prime} \in \F-\{0\}.
$$
Choose $a^\prime \in \mathbb{F} -\{0\}$ so as to satisfy the following conditions:
\begin{equation}\label{condiciones en a}
\sigma^{\prime}  \neq 0,\quad \lambda^{\prime} \neq 0,\quad \mbox{and} \quad 2{\sigma^{\prime}}^2-3\zeta^{\prime}\lambda^{\prime}\neq 0.
\end{equation}
Now, choose some $a\in\Bbb F-\{0\}$ and compute $L_a(D^\prime)$ to obtain,
$$
\aligned
L_a(D^{\prime})
&=(\,\zeta^{\prime}-2a \sigma^{\prime}+a^2 \lambda^{\prime},a(3\zeta^{\prime}-3a \sigma^{\prime}+a^2\lambda^{\prime}),
\\
&\quad\qquad
\sigma^{\prime}-a \lambda^{\prime},\lambda^{\prime},
a^2(6\zeta^{\prime}-4a \sigma^{\prime}+a^2\lambda^{\prime})\,).
\endaligned
$$
The new components of $L_{a}(D^{\prime})=(\zeta^{\prime \prime},\eta^{\prime \prime},\sigma^{\prime \prime},\lambda^{\prime \prime},\mu^{\prime \prime})$:
$$
\zeta^{\prime \prime}\!=\!\zeta^{\prime}-2a \sigma^{\prime}+a^2 \lambda^{\prime},\quad \eta^{\prime \prime}\!=\!a(3\zeta^{\prime}-3a \sigma^{\prime}+a^2\lambda^{\prime}),
$$
$$
\mu^{\prime \prime}\!=\!a^2(6\zeta^{\prime}-4a \sigma^{\prime}+a^2\lambda^{\prime}),
$$
are polynomials in $a$. Define the quadratic polynomials,
$$
p_1=\zeta^{\prime}-2\sigma^{\prime}x+\lambda^{\prime} x^2,\quad p_2=3\zeta^{\prime}-3\sigma^{\prime}x+\lambda^{\prime} x^2,
$$
$$
p_3=6 \zeta^{\prime}-4\sigma^{\prime} x+\lambda^{\prime}x^2.
$$
Then, $\zeta^{\prime \prime}=p_1(a)$, $\eta^{\prime \prime}=a\,p_2(a)$ 
and $\mu^{\prime \prime}=a^2\,p_3(a)$. These polynomials are relatively prime, 
since $2\,p_2-p_1-p_3=-\zeta \in \mathbb{F} -\{0\}$. 
We claim that there exists $a \in \mathbb{F} -\{0\}$ such that $p_1(a)p_2(a) \neq 0=p_3(a)$, 
in which case, $\zeta^{\prime \prime} \eta^{\prime \prime} \neq 0
= \mu^{\prime \prime}$.
Observe that $p_3-p_2=-\sigma^{\prime}x+3\zeta^{\prime}$. 
If any zero of $p_3$ were a zero of $p_2$, then 
$0=p_2\left(3\zeta^{\prime}/\sigma^{\prime}\right)
=p_3\left(3\zeta^{\prime}/\sigma^{\prime}\right)=
-6\zeta^{\prime}+9\,{\zeta^{\prime}}^2\, \lambda^\prime/{\sigma^{\prime}}^2$. 
Whence, $2{\sigma^{\prime}}^2=3\zeta^{\prime}\lambda^{\prime}$, contradicting the hypothesis that $2{\sigma^{\prime}}^2-3\zeta^{\prime}\lambda^{\prime} \neq 0 $, given in \eqref{condiciones en a}. Therefore, there exists $a \in \mathbb{F} -\{0\}$ such that $p_2(a) \neq 0$ and $p_3(a)=0$. If $p_1(a)=0$, we are done, because 
$p_1(a)=\zeta^{\prime \prime}=0$ and $D$ becomes equivalent to 
$(0,\eta^{\prime \prime},\sigma^{\prime \prime},\lambda^{\prime \prime},0)$ as claimed. 
Otherwise, $p_1(a) \neq 0$ implies that $\zeta^{\prime} \neq 0$ and $D$ is equivalent to $
D^{\prime \prime}=
(\zeta^{\prime \prime},\eta^{\prime \prime},\sigma^{\prime \prime},\lambda^{\prime \prime},0)$, with $\zeta^{\prime \prime}\eta^{\prime \prime} \neq 0$.
Now compute
$K_{-\zeta^{\prime \prime}/(2\eta^{\prime \prime})}(D^{\prime\prime})$
to see that $D^{\prime\prime}$ (hence $D$) is equivalent to
$(0,\eta^{\prime \prime \prime},\sigma^{\prime \prime},\lambda^{\prime \prime \prime},0)$. 

\smallskip

In summary, we have proved that any $D \in \operatorname{Der}_{(-1,1,1)}(\mathfrak{sl}_2)$, 
of the form $D=(\zeta,0,\sigma,\lambda,0)$
with $\zeta \neq 0$ and either $\eta \neq 0$ or $\lambda \neq 0$, is equivalent to
$T^{\prime\prime\prime}=(0,\eta^{\prime \prime \prime},\sigma^{\prime \prime \prime},\lambda^{\prime \prime \prime},0)$; that is,
$$
D\sim\begin{pmatrix}
0 & \eta^{\prime \prime \prime} & \sigma^{\prime \prime \prime} \\
2\sigma^{\prime \prime \prime} & 0 & \lambda^{\prime \prime \prime} \\
2 \eta^{\prime \prime \prime} & 0 & 0
\end{pmatrix}\!.
$$
Now suppose $\sigma=\lambda=0$, then $D=(\zeta,0,0,0,0)$. Applying $J_{-1,c}$ to $D$, we get:,
$$
J_{-1,c}(\zeta,0,0,0,0)=
(-\zeta,0,6c\zeta,24c^{2}\zeta,0),
$$
where $-\zeta\neq 0$, $6c\zeta \neq 0$ and $24c^2\zeta\neq 0$. By the above argument, we may deduce again that $(-\zeta,0,6c\zeta,24c^{2}\zeta,0)$ is equivalent to $(0,\eta^{\prime \prime},\sigma^{\prime \prime},\lambda^{\prime \prime},0)$.  

Finally, we apply the diagonal automorphism: $H \mapsto H$, $E \mapsto \nu E$, $F \mapsto \nu^{-1} F$, to  $(0,\eta^{\prime \prime},\sigma^{\prime \prime},\lambda^{\prime \prime},0)$, in order to obtain our result. 
\end{proof}
\noindent Therefore, from now on, we consider any generalized derivation $D$ as
$$
D=(0,\eta,\sigma,\lambda,0)=\begin{pmatrix}
0 & \eta & \sigma \\
2\sigma & 0 & \lambda \\
2 \eta & 0 & 0
\end{pmatrix}\!,\quad \eta,\sigma,\lambda \in \mathbb{F}.
$$
\begin{theorem}\label{theorem clasificacion}{\sl
Let $\F$ be an algebraically closed field of characteristic zero. Let $D \in \Der_{(-1,1,1)}(\mathfrak{sl}_2(\F))$ be a generalized derivation and let $T \in \End_{\mathbb{F}}(\mathfrak{sl}_2(\F)[D])$ be a linear map defined by $T(x+\zeta D)=x-\zeta D$ for all $x \in \mathfrak{sl}_2(\F)$ and for all $\zeta \in \mathbb{F}$. Let $\mathfrak{sl}_2(\F)[D]=\mathfrak{sl}_2(\F) \oplus \F D$ be the Hom-Lie algebra over $\F$, defined by $D$.

(i)\, If $\operatorname{rank}(D)\!=\!1$, then $\mathfrak{sl}_2(\F) \!\oplus \!\mathbb{F}D \!\simeq \! \mathfrak{sl}_2(\F) \! \oplus \!\mathbb{F} \!\!\begin{pmatrix} 0 & 0 & 0 \\ 0 & 0 & 1 \\ 0 & 0 & 0 \end{pmatrix}$.

(ii)\, If  $\operatorname{rank}(D)=2$, then
\begin{itemize}

\item $\mathfrak{sl}_2(\F)[D] \simeq \mathfrak{sl}_2(\F) \oplus \mathbb{F} \begin{pmatrix} 0 & 0 & 1 \\ 0 & 0 & 0 \\ 2 & 0 & 0 \end{pmatrix}$, or

\item $\mathfrak{sl}_2(\F)[D] \simeq \mathfrak{sl}_2(\F) \oplus \mathbb{F} \begin{pmatrix} 0 & 1 & \sigma \\ 2\sigma & 0 & 0 \\ 2 & 0 & 0 \end{pmatrix}$, with $\sigma \neq 0$.

\end{itemize}

(iii) If $\operatorname{rank}(D)=3$, then

\begin{itemize}

\item $\mathfrak{sl}_2(\F)[D] \simeq \mathfrak{sl}_2(\F) \oplus \mathbb{F}\begin{pmatrix} 0 & 1 & 0 \\ 0 & 0 & \lambda \\ 2 & 0 & 0 \end{pmatrix}$, with $\lambda \neq 0$, or

\item $\mathfrak{sl}_2(\F)[D] \simeq \mathfrak{sl}_2(\F) \oplus \mathbb{F}\begin{pmatrix} 0 & 1 & \sigma \\ 2\sigma & 0 & \lambda \\ 2 & 0 & 0 \end{pmatrix}$, with $\sigma \lambda \neq 0$.

\end{itemize}
}

\end{theorem}

\begin{proof}

{\bf Case 1.} 
{\boldmath $\operatorname{rank}(D)=1$.}
Suppose $\dim_{\mathbb{F}}\operatorname{Span}\{v_1,v_2,v_3\}=1$. Under this assumption, it is not difficult to prove that any $D \in \operatorname{Der}_{(-1,1,1)}(\mathfrak{sl}_2)$ 
with $\operatorname{rank}(D)=1$, is equivalent to $\left(\begin{smallmatrix}
\,0 & \,\,\, 0 & 0 \\
\,0 & \,\,\, 0 & \,\,\lambda \\
\,0 & \,\,\, 0 & 0
\end{smallmatrix}\right)$ with $\lambda \neq 0$. Applying a diagonal automorphism to $(0,0,0,\lambda,0)$, we get that any $D \in \operatorname{Der}_{(-1,1,1)}(\mathfrak{sl}_2)$ 
with $\operatorname{rank}(D)=1$, is equivalent to $(0,0,0,1,0)$.

\medskip
{\bf Case 2.} 
{\boldmath $\operatorname{rank}(D)=2$.}
Notice that, $v_1\wedge v_2\wedge v_3 =2\,{\eta}^2\,\lambda\,H\wedge E\wedge F$.
Thus, $\det(D)=2\,{\eta}^2\,\lambda$. 
Then, either $\eta=0$ or $\lambda=0$, under the assumption $\operatorname{rank}(D)=2$.
If $\eta=0$, then $v_2\wedge v_3=v_1\wedge v_2=0$, and 
$v_3\wedge v_1 = 2\,{\sigma}^2\,H\wedge E$. In particular, $\sigma\ne 0$,
otherwise, $\operatorname{rank}(D)=1$. On the other hand, 
the condition $\lambda=0$ under the assumption $\operatorname{rank}(D)=2$,
requires that either $\eta\ne 0$ or $\sigma\ne 0$. That is, $(\eta,\sigma)\ne (0,0)$.
Thus, we obtain two types of possible canonical forms: either $T$ is equivalent to,
\begin{equation}\label{rango 2 sl_2}
\begin{pmatrix}
 0 & 0 & \sigma \\
2\sigma & 0 & \lambda \\
0 & 0 & 0
\end{pmatrix}\!,\ \sigma \neq 0,
\mbox{ or to }
\begin{pmatrix}
  0 & \eta & \sigma \\
2\sigma & 0 & 0 \\
2\eta & 0 & 0
\end{pmatrix}\!,  
\begin{pmatrix}
\eta\\ \sigma
\end{pmatrix} \neq 0.
\end{equation}
Now, the characteristic polynomial of $D$ corresponding to the first type above
is $\chi_D=x^3$, whereas the characteristic polynomial of $T$ corresponding to the
second type is $\chi_D=x(x^2-4\,\eta\sigma)$. Clearly, these two polynomials coincide
if $\eta\sigma=0$. We shall now prove that if this is the case, then
both types given in \eqref{rango 2 sl_2}, can be made equivalent to $(0,0,\sigma,0,0)$, with $\sigma \neq 0$.
\smallskip

First notice that $D=(0,\eta,0,0,0)$ is equivalent to $D^\prime=(0,0,\sigma^\prime,0,0)$, whenever $\eta\ne 0$ and $\sigma^{\prime}\ne 0$. Indeed, take $c\in\Bbb F - \{0\}$ and see that $J_{1,c}(D)=(0,0,c^2\eta,8c^3\eta,0)$. On the other hand, let $a \in \F - \{0\}$, then $K_a(D^{\prime})=(0,0,\sigma^{\prime},-4a\sigma^{\prime})$. The elections $a=-2\sqrt{\sigma^{\prime}}$, and $c=\sqrt{\frac{\sigma^{\prime}}{\eta}}$, yields $J_{1,c}(D)=K_a(D^{\prime})$, which proves our assertion. This also proves that $D^\prime=(0,0,\sigma^\prime,0,0)$ is equivalent to 
$D^{\prime\prime}=(0,0,\sigma^{\prime\prime},\lambda^{\prime\prime},0)$ whenever $\sigma^\prime\ne 0$. Now, applying a diagonal automorphism, we obtain that $(0,0,\sigma,0,0)$, with $\sigma \neq 0$, is equivalent to $(0,0,1,0,0)$. It remains the case when $\eta\sigma \neq 0$. For that, we apply a diagonal automorphism to $D=(0,\eta,\sigma,0,0)$, with $\eta\sigma \neq 0$, in order to obtain that such $D$ is equivalent to $(0,1,\sigma^{\prime},0,0)$, with $\sigma^{\prime} \neq 0$.
Therefore,
we end up with two different canonical forms; namely,
that corresponding to the data $D=(0,0,1,0,0)$,
and that corresponding to the data $D=(0,1,\sigma,0,0)$
with $\sigma\ne 0$.

\medskip

{\bf Case 3.} 
{\boldmath $\operatorname{rank}(D)=3$.}
In this case, $D$ is invertible, and $\det(D)=2\eta^2 \lambda \neq 0$. 
The characteristic polynomial of $D$ is $x^3-4\eta \sigma x-2 \eta^2\lambda$.
If $4 \eta \sigma=0$, then $\sigma=0$ and $D=(0,\eta,0,\lambda,0)$,
with characteristic polynomial $x^3-2\eta^2\lambda$. On the other hand,
if $\sigma \neq 0$, then $D=(0,\eta,\sigma,\lambda,0)$ 
and its characteristic polynomial is 
$x^3-4 \eta \sigma x-2\lambda$, with $\eta \sigma\lambda \neq 0$.
Therefore, we obtain
two non-isomorphic canonical forms for $D$; either,
$$
\begin{pmatrix}
0 & 1 & 0 \\
0 & 0 & \lambda \\
2 & 0 & 0
\end{pmatrix}\!\!,\,\,\sigma \neq 0,
\ \ \text{or}\ \ 
\begin{pmatrix}
0 & 1 & \sigma \\
2 \sigma & 0 & \lambda \\
2  & 0 & 0
\end{pmatrix}\!\!,\,\,\sigma \lambda \neq 0.
$$
\end{proof}

\subsection{Representation theory for the Hom-Lie algebra $\mathfrak{sl}_2(\F)[D]$}

The known definition of a representation of a Hom-Lie algebra
$(\mathfrak{g},[\cdot,\cdot]_{\mathfrak{g}},T)$ in a vector space $V$
can be found in \cite{Benayadi}. 
That definition is slightly different from what one is used to find in classical representation theory.

\begin{Def}\label{Def Hom-Lie representacion}{\rm
Let $(\mathfrak{g},[\cdot,\cdot]_{\mathfrak{g}},T)$ be a Hom-Lie algebra and let $V$ be a vector space. 
A representation of  $(\mathfrak{g},[\cdot,\cdot]_{\mathfrak{g}},T)$ on $V$ with respect to 
$L \in \mathfrak{gl}(V)$, is a linear map $\rho:\mathfrak{g} \to \mathfrak{gl}(V)$, satisfying,
$$
\rho([x,y]_{\mathfrak{g}}) \circ L=\rho(T(x)) \rho(y)-\rho(T(y)) \rho(x),\quad \forall x,y \in \mathfrak{g}.
$$
}
\end{Def}
\noindent
One may then prove the following result which basically states that a
Hom-Lie algebra representation $\rho:\g \to \gl(V)$ with respect to $L:V\to V$
is equivalent to the existence of a Hom-Lie algebra structure on $\g \oplus V$, where
the product
$[\,\cdot\,,\,\cdot\,]_{\g \oplus V}$ is defined 
by $[x+u,y+v]_{\g \oplus V}=[x,y]_{\g}+\rho(x)(v)-\rho(y)(u)$
and $S=T\oplus L$
satisfies the  Hom-Lie Jacobi identity. In this way,
$\g$ becomes a Hom-Lie subalgebra and $V$ a Hom-Lie ideal.
Now, see \cite{Benayadi} for the proof of the following:

\begin{Prop}{\sl
Let $(\g,[\cdot,\cdot]_{\g},T)$ be a Hom-Lie algebra and let $V$ be a vector space. Let $\rho:\g \rightarrow \mathfrak{gl}(V)$ be a representation of this Hom-Lie algebra with respect to $L \in \mathfrak{gl}(V)$
in the sense of {\bf Definition \ref{Def Hom-Lie representacion}}.
Let $S \in \End_{\mathbb{F}}(\g \oplus V)$, be the linear map defined by $S(x+v)=T(x)+L(v)$, for all $x \in \g$ and for all $v \in V$. Let $[\,\cdot\,,\,\cdot\,]_{\g \oplus V}:(\g \oplus V) \times (\g \oplus V) \to \g \oplus V$, be the  product defined by $[x+u,y+v]_{\g \oplus V}=[x,y]_{\g}+\rho(x)(v)-\rho(y)(u)$, for all $x,y \in \g$ and for all $u,v \in V$. Then $(\g \oplus V,[\,\cdot\,,\,\cdot\,]_{\g \oplus V},S)$ is a Hom-Lie algebra.}
\end{Prop}

We shall now address one of the main goals of this work.
In the following result we characterize the Hom-Lie modules of 
the Hom-Lie algebras of the form $\mathfrak{sl}_2(\F)[D]$, 
with $D \in \Der_{(-1,1,1)}(\mathfrak{sl}_2(\F))$. 
Basically, our result allows us to conclude that
any non-trivial irreducible representation of $\mathfrak{sl}_2(\F)[D]$
in a vector space $V$ ({\it ie\/,} one for which the image of $D$ is non-zero)
must only have $\dim V=3$.

\begin{theorem}\label{theorem-rep}{\sl
Let $\mathfrak{sl}_2(\F)=\operatorname{Span}_{\F}\{H,E,F\}$ be the rank-one simple Lie algebra
over an algebraically closed field $\F$ of characteristic zero and let $D=(\zeta,\eta,\sigma,\lambda,\mu) \in \Der_{(-1,1,1)}(\mathfrak{sl}_2(\F))$ be a generalized derivation. Let $T \in \End_{\mathbb{F}}(\mathfrak{sl}_2(\F))$ be the linear map defined by $T(x+\xi D)=x-\xi D$, for all $x \in \mathfrak{sl}_2(\F)$ and for all $\xi \in \mathbb{F}$. Let $\mathfrak{sl}_2(\F)[D]$ be the Hom-Lie algebra defined by $D$. Let $V$ be a finite dimensional vector space over $\F$ and let $\rho:\mathfrak{sl}_2(\F)[D] \to \mathfrak{gl}(V)$ be a non-zero Hom-Lie algebra representation with respect to the identity map $\operatorname{Id}_V$. Then,

\begin{enumerate}

\item The representation $\rho:\mathfrak{sl}_2(\F)[D] \to \gl(V)$ is completely reducible.

\item If $\rho:\mathfrak{sl}_2(\F)[D] \to \gl(V)$ is irreducible, then 
$\rho|_{\mathfrak{sl}_2(\F)}:\mathfrak{sl}_2(\F) \to \mathfrak{gl}(V)$ is an irreducible
representation of the Lie algebra $\mathfrak{sl}_2(\F)$.

\item Let $\rho:\mathfrak{sl}_2(\F)[D] \to \gl(V)$ be irreducible. 
If $\rho(D) \neq 0$, then $\dim_{\mathbb{F}}V=3$. In this case, there is a basis $\{v_0,v_1,v_2\}$ of $V$
such that $\rho(H)(v_0)=2v_0$, $\rho(H)(v_1)=0$ and $\rho(H)(v_2)=-2v_2$. Furthermore, the matrix of $\rho(D)$ in the basis $\{v_0,v_1,v_2\}$, has the block form,
$$
\rho(D)=
\begin{pmatrix}
-\zeta & -2 \sigma & -\lambda\\
-\eta & 2 \zeta & \sigma\\
-\mu & 2 \eta & -\zeta
\end{pmatrix}\!.
$$
If $\rho(D)=0$, then $D=0$.
\end{enumerate}
}
\end{theorem}

\begin{proof}
$\,$
\smallskip

{\boldmath $(1)$} Let $\rho:\mathfrak{sl}_2(\F)[D] \to \mathfrak{gl}(V)$ be as in the statement.  Then,
$$
\rho([x,y])=\rho(T(x))\rho(y)-\rho(T(y))\rho(x),\quad \forall x,y \in \mathfrak{sl}_2(\F)[D].
$$ 
In particular, $\rho(D(x))=-\rho(D)\rho(x)-\rho(x)\rho(D)$, for all $x \in \mathfrak{sl}_2(\F)$. Note that $\rho|_{\mathfrak{sl}_2(\F)}$ is a representation of the Lie algebra $\mathfrak{sl}_2(\F)$ in $V$. Let $U$ be an $\mathfrak{sl}_2(\F)[D]$-submodule of $V$. Since $\rho|_{\mathfrak{sl}_2(\F)}$ is 
completely reducible, there exists an $\mathfrak{sl}_2(\F)$-submodule $W \subset V$, such that $V=U \oplus W$. We claim that $\rho(D)(W) \subset W$.
\smallskip

There are linear maps $T \in \operatorname{Hom}_{\mathbb{F}}(W,U)$ and $S \in \operatorname{End}_{\mathbb{F}}(W)$, such that $\rho(D)(w)=T(w)+S(w)$, for all $w \in W$. Let $x \in \mathfrak{sl}_2(\F)$ and $w \in W$. Then,
$$
\aligned
\rho(x)(T(w))+\rho(x)(S(w))&=\rho(x)\circ \rho(D)(w)\\
\,&=-\rho(D(x))(w)-\rho(D)\circ \rho(x)(w)\\
\,&=-T(\rho(x)(w))-\rho(D(x))(w)-S(\rho(x)(w)).
\endaligned
$$
It follows that $\rho(x)(T(w))=-T(\rho(x)(w))$, for all $w \in W$ and all $x \in \mathfrak{sl}_2(\F)$. Then $T \in \operatorname{Hom}_{\mathbb{F}}(W,U)$ is a linear map satisfying $\rho(x) \circ T=-T \circ \rho(x)$, for all $x \in \mathfrak{sl}_2(\F)$. This also implies that $\Ker T$ and $\Im T$ are $\mathfrak{sl}_2(\F)$-submodules of $W$ and $U$, respectively. Thus, there are $\mathfrak{sl}_2(\F)$-submodules, $W' \subset W$ and $U' \subset U$, such that $W=\Ker T \oplus W'$ and $U=\Im T \oplus U'$. Therefore, $T|_{W'}:W' \to \Im T$ is an isomorphism of $\mathfrak{sl}_2(\F)$-modules. If $W'=\{0\}$, then $W=\Ker T$ and $\rho(D)(w)=S(w) \in W$, for all $w \in W$. 

\smallskip

So, assume $W' \neq \{0\}$. With no loss of generality, we may assume that
$T \in \operatorname{Hom}_{\F}(W,U)$ is an isomorphism.
Let $W'' \subset W$ be an irreducible $\mathfrak{sl}_2(\F)$-submodule, then $T(W'')$ is an irreducible $\mathfrak{sl}_2(\F)$-submodule of $U$. Assume that both $W$ and $U$ are irreducible $\mathfrak{sl}_2(\F)$-modules. We write $W=\operatorname{Span}_{\mathbb{F}}\{w_0,w_1,\ldots,w_m\}$ (resp. $U=\operatorname{Span}_{\mathbb{F}}\{u_0,u_1,\ldots,u_m\}$), where:
$$
\aligned
\rho(H)(w_k) & = (m-2k)w_k,
\\
\rho(E)(w_k) & = (m-k+1)w_{k-1},
\\
\rho(F)(w_k) & =(k+1)w_{k+1},
\endaligned
$$
(see \cite{Hum}, {\bf Lemma 7.2}) for all $0 \leq k \leq m$, with the convention, $w_{-1}=w_{m+\ell}=0$, for $\ell \geq 1$. 
Analogous expressions are satisfied in $U$ for the $u_k$'s. Now,
$$
T \circ \rho(H)(w_k)=-\rho(H) \circ T(w_k)\ \Rightarrow\ \rho(H)(T(w_k))=-(m-2k)T(w_k).
$$
This means that the weight of $T(w_k)$ is $-(m-2k)=m-2(m-k)$. Thus, there exists a scalar $\xi_k \in \mathbb{F}$, such that $T(w_k)=\xi_k u_{m-k}$, for all $0 \le k \le m$. For $m-k$, we get $T(w_{m-k})=\xi_{m-k}u_k$.
Similarly, from $T \circ \rho(E)(w_k)=-\rho(E) \circ T(w_k)$, we obtain,
\begin{equation*}\label{rep1}
(m-k+1)\xi_{k-1}u_{m-k+1}=-(k+1)\xi_k u_{m-k-1}.
\end{equation*}
And from  $T \circ \rho(F)(w_k)=-\rho(F) \circ T(w_k)$, we get,
\begin{equation*}
(k+1)\xi_{k+1}u_{m-k-1}=-(m-k+1)\xi_k u_{m-k+1}.
\end{equation*}
Since $m-k-1 \neq m-k+1$, we conclude that $u_{m-k-1}$ and $u_{m-k+1}$ are linearly independent.
Thus,
$$
\aligned
(m-k+1)\xi_{k-1}&=-(k+1)\xi_{k}=0,
\\
-(m-k+1)\xi_{k}&=(k+1)\xi_{k+1}=0, \qquad 0 \leq k \leq m.
\endaligned
$$
Whence $\xi_k=0$, for all $k$ and therefore, $T=0$. This also shows that $\rho(D)(W) \subset W$ and therefore $W$ is a $\mathfrak{sl}_2(\F)[D]$-module complementary to $U$ in $V$. This argument shows that $\rho$ is completely reducible.
\medskip

{\boldmath $(2)$} Now, let $\rho:\mathfrak{sl}_2(\F)[D] \to \gl(V)$ be irreducible as a representation of the given Hom-Lie algebra $\mathfrak{sl}_2(\F)[D]$ in $V$. In particular, $V$ is an $\mathfrak{sl}_2(\F)$-module. Let $\{0\} \neq U \subset V$ be an $\mathfrak{sl}_2(\F)$-submodule. Then, there there exists an $\mathfrak{sl}_2(\F)$-submodule $W \subset V$, such that $V=U \oplus W$. We then obtain linear maps 
$S' \in \operatorname{End}_{\mathbb{F}}(U)$ and 
$T' \in \operatorname{Hom}_{\mathbb{F}}(U,W)$, such that $\rho(D)(u)=S'(u)+T'(u)$, for all $u \in U$. From $\rho(D(x))=-\rho(D)\circ \rho(x)-\rho(x) \circ \rho(D)$, for all $x \in \mathfrak{sl}_2(\F)$, it follows $\rho(x) \circ T^{\prime}=-T^{\prime} \circ \rho(x)$, for all $x \in \mathfrak{sl}_2(\F)$, where $T^{\prime}:U \to W$. Using the same argument as in \textbf{(1)} above, one concludes that $T'=0$ and $\rho(D)(U) \subset U$. From the beginning, we made the assumption that $V$ is an irreducible $\mathfrak{sl}_2(\F)[D]$-module, therefore $V=U$ and $W=\{0\}$.
\medskip

{\boldmath $(3)$} It follows that if $\rho:\mathfrak{sl}_2(\F)[D] \to \gl(V)$ is irreducible
in the Hom-Lie category (with respect to $\operatorname{Id}_{V}$), then $V$ is an irreducible $\mathfrak{sl}_2(\F)$-module. Thus, there is a basis $\{v_0,v_1,...,v_m\}$ of $V$, for which,
$\rho(H)(v_k)=(m-2k)v_k$, $\rho(E)(v_k)=(m-k+1)v_{k-1}$, and $\rho(F)(v_k)=(k+1)v_{k+1}$, for all $0 \leq k \leq m$.
It only remains to determine how $\rho(D)$ acts in the basis $\{v_0,v_1,\ldots,v_m\}$. Let $a_{i,j}$ be the scalars for which 
$$
\rho(D)(v_j)=\overset{m}{\underset{i=0}\sum}a_{i,j}v_i,\quad 0\le j \le m.
$$
We shall now determine the $(m+1) \times (m+1)$ matrix $A \in \operatorname{Mat}_{m+1 \times m+1}(\mathbb{F})$, whose $i,j$ entry is $a_{i,j}$.
Let $D=(\zeta,\eta,\sigma,\lambda,\mu)$ be a $(-1,1,1)$-derivation of $\mathfrak{sl}_2(\F)$. Then,
$$
\aligned
D(H) & = 2 \zeta H+2 \sigma E+2 \eta F,
\\
D(E) & = \eta H-\zeta E+\mu F,
\\
D(F) & = \sigma H+\lambda E-\zeta F.
\endaligned
$$
From $\rho(H)\circ \rho(D)(v_j)=-\rho(D(H))(v_j)-\rho(D)\rho(H)(v_j)$, we obtain,
\begin{equation}\label{repn1}
\aligned
\sum_{i=0}^m 
& 
\left(\,m-(i+j)\,\right)\,a_{i,j}\,v_i 
+(m-2j)\zeta v_j +\\
&
\quad
+(m-j+1)\sigma \, v_{j-1}+(j+1)\eta \, v_{j+1}=0, \quad 0 \le j \le m.
\endaligned
\end{equation}
We deduce the following:
\begin{eqnarray}\label{repn2}
i=j \Rightarrow &(m-2j)a_{j,j}&=-(m-2j)\zeta,\,\,\,\forall j,\\
\label{repn3}
i=j+1 \Rightarrow & (m-(2j+1))a_{j+1,j}&=-(j+1)\eta,\,\,\,\forall j,\\
\label{repn4}
i=j-1 \Rightarrow & (m-(2j-1))a_{j-1,j}&=-(m-j+1)\sigma,\,\,\,\forall j,
\end{eqnarray}
Similarly, from $\rho(E)\circ \rho(D)(v_j)=-\rho(D(E))(v_j)-\rho(D) \circ \rho(E)(v_j)$, 
we get for $0\le j\le m$,
\begin{equation}\label{repn5}
\aligned
\sum_{i=0}^{m-1}\left(\,(m-i)a_{i+1,j}+(m-j+1)a_{i,j-1}\,\right)v_i+(m-2j)\eta\, v_j&\,\\
-(m-j+1)\zeta \, v_{j-1}+(j+1) \mu \, v_{j+1}+(m-j+1)a_{m,j-1}\,v_m&=0.
\endaligned
\end{equation}
For $i=j-1$ in \eqref{repn5}, we obtain,
$$
((m-j+1)a_{j,j}+(m-j+1)a_{j-1,j-1})v_{j-1}-(m-j+1)\zeta\,v_{j-1}=0,
$$
and therefore,
\begin{equation}\label{repn7}
a_{j,j}+a_{j-1,j-1}=\zeta,\quad 1 \leq j \leq m.
\end{equation}
Now observe that \eqref{repn7} yields, $a_{j-1,j-1}=a_{j+1,j+1}$, for all $j$. 
If   in \eqref{repn2} we make the changes 
$j \to j-1$ and $j \to j+1$, we get,
\begin{eqnarray}
(m-(j-1))a_{j-1,j-1}&=-(m-(j-1))\zeta,\label{f1}\\
(m-(j+1))a_{j+1,j+1}&=-(m-(j+1))\zeta\label{f2}.
\end{eqnarray}
By subtracting these two expressions, we get $a_{j-1,j-1}=-\zeta$, and therefore \eqref{repn7} implies $a_{j,j}=2\zeta$. This also proves that $a_{j+1,j+1}=a_{j-1,j-1}=-\zeta$, if and only if $a_{j,j}=2 \zeta$, for $1 \leq j \leq m$. On the other hand, if $m$ is odd, it follows from \eqref{repn2} that,
$a_{j,j}=-\zeta$ for all $j$, which implies $\zeta=0$. Thus, let us assume that $m=2\ell$, with $\ell \in \mathbb{N}$. From \eqref{repn2}, we know that $(m-2j)a_{j,j}=-(m-2j)\zeta$, for all $j$. Whence, $a_{j,j}=-\zeta$ for all $j \neq \ell$, and thus $a_{\ell,\ell}=2\zeta$. If $\ell-2 \geq 0$ then $2\zeta=a_{\ell,\ell}=a_{\ell-2,\ell-2}=-\zeta$, so, $\zeta=0$. This proves that if $m\neq 2$, then $\zeta=0$.
\smallskip

\noindent In addition, for $j=m$ in \eqref{repn5}, it follows that $a_{m,m-1}=m\,\eta$. On the other hand, applying \eqref{repn3} with $j=m-1$, we get $(m-1)a_{m,m-1}=m\,\eta$. Therefore, $(m-2)\eta=0$.
\smallskip

\noindent Now, from $\rho(F) \circ \rho(D)(v_j)=-\rho(D(F))(v_j)-\rho(D)\rho(F)(v_j)$, we obtain,
\begin{equation}\label{repn8}
\aligned
(j+1) 
&
a_{0,j+1}v_0  
+\sum_{i=1}^m(ia_{i-1,j}+(j+1)a_{i,j+1})v_i
\\
& +
(m-2j)\sigma v_j+(m-j+1)\lambda v_{j-1}-(j+1)\zeta v_{j+1}=0.
\endaligned
\end{equation}
Taking $j=0$, we get $a_{0,1}+m\,\sigma=0$. Combining this with \eqref{repn4} and $j=1$, we get $(m-2)\,\sigma=0$.

\smallskip
\noindent 
Let us assume $m \neq 2$, then $\zeta=\eta=\sigma=0$. Then, both \eqref{repn5} and \eqref{repn8}, respectively, take the form:
\begin{multline}\label{repn5a}
\sum_{i=0}^{m-1}((m-i)a_{i+1,j}+(m-j+1)\,a_{i,j-1})v_{i}\\
+(j+1)\mu v_{j+1}+(m-j+1)\,a_{m,j-1}v_m=0
\end{multline}
\begin{equation}\label{repn8a}
(j+1)a_{0,j+1}v_0+\sum_{i=1}^m(ia_{i-1,j}+(j+1)a_{i,j+1})v_{i}+(m-j+1)\lambda v_{j-1}=0,
\end{equation}
for all $j$. Since $m \neq 2$, from \eqref{repn5a}, it follows that for $j=1$ we obtain
$$
\sum_{i=0}^{m-1}((m-i)a_{i+1,1}+m\,a_{i,0})v_{i}+2\mu v_{2}+m\,a_{m,0}v_m=0,
$$
then $a_{m,0}=0$. Similarly, for $j=m-1$ on \eqref{repn8a}, we have 
$$
m\,a_{0,m}v_0+\sum_{i=1}^m(ia_{i-1,m-1}+m\,a_{i,m})v_{i}+2\lambda v_{m-2}=0,
$$
then $a_{0,m}=0$. On the other hand, from \eqref{repn5a} and \eqref{repn8a}, we obtain, 
\begin{eqnarray}
\label{repn9a}
(m-i)a_{i+1,j}+(m-j+1)a_{i,j-1}&=0,&\,i \notin\{j+1,m\},\\
\label{repn9b}
ia_{i-1,j}+(j+1)a_{i,j+1}&=0,&\,i \notin\{0,j-1\}.
\end{eqnarray}
Letting $j \to j+1$ in \eqref{repn9a}  and $j \to j-1$ in \eqref{repn9b}, respectively, we get:
\begin{eqnarray}
\label{repn10a}
(m-i)a_{i+1,j+1}+(m-j)a_{i,j}&=0,&\,i \notin\{j+2,m\},\\
\label{repn10b}
ia_{i-1,j-1}+ja_{i,j}&=0,&\,i \notin\{0,j-2\}.
\end{eqnarray}
On the other hand, since $\zeta=\eta=\sigma=0$, it follows from \eqref{repn1} that, 
$$
\sum_{i=0}^m\left(\,m-(i+j)\,\right)a_{i,j}v_i=0,\quad\,0 \leq j \leq m.
$$
Then, $a_{i,j}=0$, for all $i,j$ satisfying $i+j \neq m$. Let $k$ and $\ell$ be such that $k+\ell=m$. 
Then, $(k+1)+(\ell+1) \neq m$ and $(k-1)+(\ell-1) \neq m$, 
so that $a_{k+1,\ell+1}=a_{k-1,\ell-1}=0$. 
We have already proved that $a_{0,m}=a_{m,0}=0$. Thus, assume $k\ne 0$ and $k\ne m$.
If $k \neq \ell \pm 2$, 
it follows from \eqref{repn10a} and \eqref{repn10b} that $(m-\ell)a_{k,\ell}=\ell a_{k,\ell}=0$.
Thus, $a_{k,\ell}=0$. If $k=\ell+2$, then $k \neq \ell-2$ and $m=k+\ell=2\ell+2$.
Since $m \neq 2$, then $\ell \neq 0$. 
It follows from \eqref{repn10b} that $a_{k,\ell}=0$. 
If $k=\ell-2$, then $k \neq \ell+2$; thus, we conclude from \eqref{repn10a} that, $0=(m-\ell)a_{k,\ell}=k a_{k,\ell}$. 
But we have assumed $k \neq 0$. Then, $a_{k,\ell}=0$.
Therefore, $a_{i,j}=0$ for all $i$ and $j$.

\smallskip
\noindent 
We have proved that if $m \neq 2$ then $\rho(D)=0$. If $m=2$, by \eqref{repn7}, \eqref{f1} and \eqref{f2}, we have $a_{0,0}=a_{2,2}=-\zeta$ and $a_{1,1}=2\zeta$, From \eqref{repn3} we obtain, $a_{1,0}=-\eta$ and $a_{2,1}=2\eta$. Similarly, it follows from \eqref{repn4} that $a_{0,1}=-2\sigma$ and $a_{1,2}=\sigma$.
On the other hand, taking $j=1$ in \eqref{repn5}, we get, $a_{2,0}=-\mu$. 
Finally, taking $j=1$ in \eqref{repn8}, we obtain $a_{0,2}=-\lambda$. We therefore conclude that
the matrix of $\rho(D)$ in the basis $\{v_0,v_1,v_2\}$, with entries $a_{i,j} \in \F$, $0 \leq i,j \leq 2$, is as claimed in the statement.
\smallskip

Finally, suppose $\rho(D)=0$. This implies that $\rho(D(x))=0$, for all $x \in \mathfrak{sl}_2(\F)$, which means that $\Im D \subset \Ker \rho|_{\mathfrak{sl}_2(\F)}$. Since $\Ker \rho|_{\mathfrak{sl}_2(\F)}$ is an ideal of $\mathfrak{sl}_2(\F)$ and $\rho \neq 0$, it follows that $\Ker \rho|_{\mathfrak{sl}_2(\F)}=\{0\}$. Therefore, $D=0$.
\end{proof}

\begin{Cor}{\sl
Let $D \neq 0$ be a $(-1,1,1)$-derivation of $\frak{sl}_2(\F)$ and consider the Hom-Lie algebra
$\frak{sl}_2(\F)\oplus\Bbb F\,D$. Any irreducible $\mathfrak{sl}_2$-Lie algebra representation 
$\rho:\mathfrak{sl}_2(\F) \to \mathfrak{gl}(V)$, can be extended to an irreducible Hom-Lie algebra representation of $\mathfrak{sl}_2(\F)[D]$ with respect to
the identity map $\operatorname{Id}_V:V\to V$,
if and only if, $\operatorname{dim}V=3$ and $\rho(D)$ is obtained as in (3)
of {\bf Thm \ref{theorem-rep}}.}
\end{Cor}

As we can see in {\bf Thm. \ref{theorem clasificacion}}, the simple Lie algebra $\mathfrak{sl}_2(\F)$ admits invertible generalized derivations of the type (-1,1,1), whereas an ordinary derivation can never be invertible in $\mathfrak{sl}_2(\F)$. On the other hand, {\bf Thm. \ref{theorem-rep}} is the analogue of Weyl's Theorem 
for Hom-Lie algebras of the form $\mathfrak{sl}_2(\F)[D]$. Thus, using {\bf Thm. \ref{theorem-rep}}, we may now state and prove the following theorem that describes the structure of the non-solvable Lie algebras that admit an invertible generalized derivation of the type $(-1,1,1)$.

\begin{theorem}\label{no solubles}{\sl
Let $\g$ be a finite dimensional non-solvable Lie algebra over an algebraically closed field $\F$ admitting an invertible generalized derivation $D \in \Der_{(-1,1,1)}(\g)$ and let $D_S \in \End_{\F}(\g)$ be its semisimple part. Then, Levi's decomposition of $\g$ is $\g=\mathfrak{sl}_2 \oplus \Rad(\g)$, with $D_S(\mathfrak{sl}_2) \subset \mathfrak{sl}_2$ and $D_S(\Rad(\g)) \subset \Rad(\g)$. Furthermore, if $[\mathfrak{sl}_2,\Rad(\g)]$ $\neq \{0\}$, and $\{v \in \Rad(\g)\,|\,[x,v]=0,\,\forall x \in \mathfrak{sl}_2\}=\{0\}$, then $\g=[\g,\g]$ and $\Rad(\g)$ decomposes into the direct sum of irreducible $\mathfrak{sl}_2$-modules, all of which correspond to the highest weight $2$.}
\end{theorem}
\begin{proof}
By {\bf Prop. 1.7} in \cite{Dorado}, it follows that $D_S \in \Der_{(-1,1,1)}$. We suppose first that $[\Rad(\g),\Rad(\g)]=\{0\}$. Let $\mathfrak{s} \subset \g$ be a complementary subspace to $\Rad(\g)$: $\g=\mathfrak{s} \oplus \Rad(\g)$. For each pair $x,y \in \mathfrak{s}$, write $[x,y]_{\mathfrak{s}} \in \s$ and $\gamma(x,y) \in \Rad(\g)$, for the components that make, $[x,y]=[x,y]_{\mathfrak{s}}+\gamma(x,y)$. On the other hand, the Lie product $[x,v] \in \Rad(\g)$, for $x \in \s$ and $v \in \Rad(\g)$, gives rise to a linear map $\rho:\mathfrak{s} \to \gl(\Rad(\g))$, through $\rho(x)(v)=[x,v]$. The Jacobi identity and the skew-symmety of the Lie bracket $[\,\cdot\,,\,\cdot\,]$, imply that $(\mathfrak{s},[\,\cdot\,,\,\cdot\,]_{\mathfrak{s}})$ is a Lie algebra, $\rho:\mathfrak{s} \to \gl(\Rad(\g))$ is a representation of $\s$ and $\gamma:\mathfrak{s} \times \mathfrak{s} \to \Rad(\g)$ is a 2-cocycle with coefficients in the representation $\rho$. Since $\mathfrak{s}$ is semisimple, then any 2-cocycle is a 2-coboundary (see \cite{Che}). Thus, there exists a linear map $\tau:\mathfrak{s} \to \Rad(\g)$, such that $\gamma(x,y)=\tau([x,y]_{\mathfrak{s}})$, for all $x,y \in \s$. Consider the semidirect sum of Lie algebras $\mathfrak{s} \ltimes_{\rho} \Rad(\g)$, with the Lie bracket $[\,\cdot\,,\,\cdot\,]_{\rho}$ given by:
$$
[x+v,y+w]_{\rho}=[x,y]_{\mathfrak{s}}+\rho(x)(w)-\rho(y)(v),\quad \forall x+v,y+w \in \mathfrak{s} \oplus \Rad(\g).
$$ 
The linear map $\g \to \mathfrak{s} \ltimes_{\rho} \Rad(\g)$, $x+v \mapsto x-\tau(x)+v$, is actually an isomorphism of Lie algebras. This proves that $\g$ is isomorphic to the semidirect sum of Lie algebras $\mathfrak{s} \ltimes_{\rho} \Rad(\g)$, where $D_{S}(\mathfrak{s}) \subset \mathfrak{s}$.
\smallskip

Now suppose that $[\Rad(\g),\Rad(\g)] \neq \{0\}$. It immediately follows that,
$$
\Rad(\g/[\Rad(\g),\Rad(\g)])=\Rad(\g)/[\Rad(\g),\Rad(\g)],
$$
is Abelian. Since $[\Rad(\g),\Rad(\g)]$ is invariant under the generalized derivation $D$, we deduce
(like in the paragraph above) that there exists a semisimple Lie subalgebra $\tilde{\mathfrak{s}} \subset \g/[\Rad(\g),\Rad(\g)]$, such that
$$
\g/[\Rad(\g),\Rad(\g)]\simeq \tilde{\mathfrak{s}} \ltimes_{\tilde{\rho}} \Rad(\g)/[\Rad(\g),\Rad(\g)],\quad \mbox{ and }\tilde{D}(\tilde{\mathfrak{s}}) \subset \tilde{\mathfrak{s}},
$$
where $\tilde{\rho}:\tilde{\mathfrak{s}} \to \gl(\g/[\Rad(\g),\Rad(\g)])$. 
\smallskip

Let $\h=\pi^{-1}(\tilde{\mathfrak{s}})$, (where $\pi:\g \to \g/[\Rad(\g),\Rad(\g)]$ is the canonical projection). Then, $\tilde{\mathfrak{s}}=\h/[\Rad(\g),\Rad(\g)]$, with $\h$ a subalgebra of $\g$. Observe that $\h$ is invariant under the generalized derivation $D$, because $\tilde{D}(\tilde{\mathfrak{s}}) \subset \tilde{\mathfrak{s}}$.  Since $\tilde{\mathfrak{s}}$ is semisimple, it follows that $\Rad(\tilde{\mathfrak{s}})=\{0\}$, thus $\Rad(\h)=[\Rad(\g),\Rad(\g)]$. Proceeding by induction on $\dim_{\F}\g$, we conclude that there exists a semisimple subalgebra $\mathfrak{s} \subset \h$, such that $\h=\mathfrak{s} \ltimes_{\rho} [\Rad(\g),\Rad(\g)]$ and $D_S(\mathfrak{s}) \subset \mathfrak{s}$. We claim that $\g=\mathfrak{s} \ltimes_{\rho} \Rad(\g)$. It is clear that $\g=\mathfrak{s}+\Rad(\g)$. Let $x \in \mathfrak{s} \cap \Rad(\g)$, then $x \in \h \cap \Rad(\g)=\Rad(\h)=[\Rad(\g),\Rad(\g)]$, from which it follows that $x \in \h \cap [\Rad(\g),\Rad(\g)]=\{0\}$. Therefore, $\g=\mathfrak{s} \ltimes_{\rho} \Rad(\g)$ and $D_S(\mathfrak{s}) \subset \mathfrak{s}$. Finally, from \cite{Burde2}, {\bf Thm. 5.12}, we have that $\mathfrak{s}=\mathfrak{sl}_2$. Thus, we may assume that the Levi-Malcev decomposition of $\g$ is $\g=\mathfrak{sl}_2(\F) \oplus \Rad(\g)$, with $D_S(\mathfrak{sl}_2) \subset \mathfrak{sl}_2$ and $D_S(\Rad(\g)) \subset \Rad(\g)$.
\smallskip

Suppose that $[\mathfrak{sl}_2,\Rad(\g)] \neq \{0\}$, which is equivalent to the fact that the representation $\rho:\mathfrak{sl}_2 \to \gl(\Rad(\g))$ is non-zero. Let us consider the Hom-Lie algebra structure on $\mathfrak{sl}_2[D_S]$, defined by $D_S$. Then, the linear map $\overline{\rho}:\mathfrak{sl}_2[D_S] \to \gl(\Rad(\g))$, defined by:
$$
\overline{\rho}(x+\xi D_S)(v)=[x,v]+\xi D_S(v),\quad \forall x \in \mathfrak{sl}_2,\quad \forall v \in \Rad(\g),\quad \forall \xi \in \F,
$$
is a Hom-Lie representation of $\mathfrak{sl}_2[D_S]$, with respect to the identity map $\Id_{\Rad(\g)}$, satisfying $\overline{\rho}(D_S)(\Rad(\g))=D_S(\Rad(\g)) \neq \{0\}$. Since $\{v \in \Rad(\g)\,|\,[x,v]=0,\,\forall x \in \mathfrak{sl}_2(\F)\}=\{0\}$,  by {\bf Thm. \ref{theorem-rep}}, it follows that $\Rad(\g)$ is a direct sum irreducible $\mathfrak{sl}_2$-modules, all of them corresponding to the highest weight $2$.  Therefore, $\Rad(\g)=[\mathfrak{sl}_2,\Rad(\g)]$, whence, $\g=[\g,\g]$ and $3$ divides $\dim_{\F}\g$.
\end{proof}

\begin{Remark}{\rm 
This result calls for the study of solvable Lie algebras
that admit invertible generalized derivations of the type $(-1,1,1)$.
It seems, however, that one requires some strong preliminary representation theory 
---such as the analogue of Weyl's Theorem that we just used
for the Hom-Lie algebra $\mathfrak{sl}_2[D_S]$---
in order to address and settle down the corresponding structure theorem
for the solvable case.
}
\end{Remark}

\begin{Remark}{\rm
The condition that $\dim_{\F}\g$ is a multiple of $3$ if $\g=\s \oplus \Rad(\g)$ is a non-solvable Lie algebra admitting an invertible generalized derivation of the type $(-1,1,1)$ with $[\s,\Rad(\g)] \neq \{0\}$, also appears for those Lie algebras admitting a \emph{periodic} generalized derivation of the type $(-1,1,1)$, (see {\bf Thm. 3.5} in \cite{Dorado}).
}
\end{Remark}

\section{Hom-Lie algebras structures on simple Lie algebras}

\noindent Let $\mathfrak{g}$ be a simple Lie algebra.
The conditions for a linear map  $T \in \End_{\mathbb{F}}(\mathfrak{g})$ to satisfy,
$$
[T(x),[y,z]]+[T(y),[z,x]]+[T(z),[x,y]]=0,\quad \forall x,y,z \in \mathfrak{g},
$$
was a question addressed and solved in \cite{Xie}. We quote the result:

\begin{theorem}[Theorem 3.3 in \cite{Xie}]{\sl
Let $\g$ be a finite-dimensional simple Lie algebra over an algebraically closed field $\F$. Let $T \in \End_{\F}(\g)$ be a linear map satisfying:
$$
[T(x),[y,z]]+[T(y),[z,x]]+[T(z),[x,y]]=0,\quad \forall x,y,z \in \mathfrak{g}.
$$ 
\begin{itemize}

\item[(i)] If $\g=\mathfrak{sl}_2(\F)$ then 
$$
T \in \operatorname{Span}_{\F}\{e_{12},e_{21},e_{33},e_{11}+e_{22},2e_{13}+e_{32},2e_{23}+e_{31}\},
$$
where $e_{ij}:\frak{sl}_2(\F)\to\frak{sl}_2(\F)$ is the linear map
defined by $e_{ij}(e_k)=\delta_{jk}e_i$
on the standard basis $\{e_1=H,e_2=E, e_3=F\}$ of $\frak{sl}_2$. Its corresponding
$3\times 3$ matrix  has a $1$ in its $(i,j)$-entry and $0$ elsewhere.

\item[(ii)] If $\g \neq \mathfrak{sl}_2(\F)$, then $T$ is a scalar multiple of the identity.

\end{itemize}
}
\end{theorem}

\begin{Remark}{\rm The proof provided in \cite{Xie} of this result, made use of the 
well known program GAP (https://www.gap-system.org).
We shall present in this section a proof from first principles of this theorem
and with no assistance of any computer program at all.}
\end{Remark}
In fact, using only the well-known theory of root systems for a simple Lie algebra, we shall prove the following:

\begin{theorem}{\sl
Let $\g$ be a finite dimensional Lie algebra over an algebraically closed field $\F$ of characteristic zero. 
Let $\operatorname{HL}(\g)$ be the vector subspace of $\End_{\F}(\g)$ defined by,
$$
\operatorname{HL}(\g)=\left\{T \in \End_{\F}(\g)\,\mid\,\sum_{\text{cyclic}}[T(x),[y,z]]=0,\mbox{ for all }x,y,z \in \g\right\}.
$$
\begin{enumerate}

\item If $\g$ is a simple Lie algebra of rank bigger than $1$, 
then $\operatorname{HL}(\g)=\F \operatorname{Id}_{\g}$. 

\item If $\g=\mathfrak{sl}_2(\F)$, then $\operatorname{HL}(\mathfrak{sl}_2(\F))$ is 
a 6-dimensional reducible $\mathfrak{sl}_2(\F)$-submodule of $\frak{gl}(\frak{sl}_2)\simeq\operatorname{End}(\frak{sl}_2)$
via the action $x. T=[\ad_{\mathfrak{sl}_2(\F)}(x),T]_{\mathfrak{gl}(\mathfrak{sl}_2(\F))}$, for all $x \in \mathfrak{sl}_2(\F)$ and $T \in {\frak{gl}}(\mathfrak{sl}_2(\F))$. In fact, 
$$
\operatorname{HL}(\mathfrak{sl}_2(\F))=\Der_{(-1,1,1)}(\mathfrak{sl}_2(\F)) \oplus \F \operatorname{Id}_{\mathfrak{sl}_2(\F)},
$$
and $\Der_{(-1,1,1)}(\mathfrak{sl}_2(\F))$ is the 
5-dimensional irreducible $\mathfrak{sl}_2(\F)$-module corresponding to the highest weight 4.

\end{enumerate}
}
\end{theorem}

\begin{proof}
{\boldmath $(1)$}  Let $\mathfrak{g}=H \oplus \underset{\alpha \in \Phi}{\bigsqcup}\mathfrak{g}_{\alpha}$ be a simple Lie algebra, where $H$ is a toral subalgebra and $\mathfrak{g}_{\alpha}=\mathbb{F}x_{\alpha}$, is the root subspace corresponding to $\alpha \in \Phi$ and $\Phi$ is a root system corresponding to $H$. Let $\pi_H:\mathfrak{g} \rightarrow H$ and $\pi_{\alpha}:\mathfrak{g} \rightarrow \mathfrak{g}_{\alpha}$, be the projections onto $H$ and $\mathfrak{g}_{\alpha}$, respectively. Let us consider the compositions,
$$
M=\pi_H \circ (T|_H):H \rightarrow  H
\ \ \text{and}\ \ N_{\alpha}=\pi_{\alpha} \circ (T|_H):H \rightarrow  \mathfrak{g}_{\alpha},
\ \ \alpha\in\Phi.
$$
We can write, $\displaystyle{T(h)=M(h)+\sum_{\alpha \in \Phi}N_{\alpha}(h)}$, for each $h \in H$. 
Similarly, consider,
$$
f_{\beta}=\pi_H \circ (T|_{\mathfrak{g}_{\beta}}):\mathfrak{g}_{\beta} \rightarrow H,
\quad\text{and}\quad
g^{\beta}_{\alpha}=\pi_{\alpha} \circ (T|_{\mathfrak{g}_{\beta}}):
\mathfrak{g}_{\beta} \rightarrow \mathfrak{g}_{\alpha},
$$
for any $\alpha$ and $\beta$ in $\Phi$.
Then $\displaystyle{T(x_{\beta})=f_{\beta}(x_{\beta})+\sum_{\alpha \in \Phi} g^{\beta}_{\alpha}(x_{\beta})}$, for all $\beta \in \Phi$. Let $h,h' \in H$ and $\alpha \in \Phi$. Then,
$$
\aligned
0& =[T(h),[h',x_{\alpha}]]+[T(h'),[x_{\alpha},h]]+[T(x_{\alpha}),[h,h']]
\\
&=
\alpha(h')[T(h),x_{\alpha}]-\alpha(h)[T(h'),x_{\alpha}]
\\
&=
\alpha(h')[M(h),x_{\alpha}]+\alpha(h')\sum_{\beta \in \Phi}[N_{\beta}(h),x_{\alpha}]
\\
&\quad -\alpha(h)[M(h'),x_{\alpha}]-\alpha(h)\sum_{\beta \in \Phi}[N_{\beta}(h'),x_{\alpha}]
\\
&=
(\alpha(h')\alpha(M(h))-M(h')\alpha(M(h)))x_{\alpha}
\\
&
\quad 
+\alpha(h')\sum_{\beta \in \Phi}[N_{\beta}(h),x_{\alpha}]
-\alpha(h)\sum_{\beta \in \Phi}[N_{\beta}(h'),x_{\alpha}].
\endaligned
$$
Observe that $[N_{\beta}(h),x_{\alpha}] \in \g_{\beta+\alpha}$, for all $\alpha,\beta \in \Phi$. Since $\alpha+\beta \neq \alpha$, for all $\beta \in \Phi$, we get from the above the following:
\begin{eqnarray}
\,\,\,\,\,\,\,\, \alpha(h') \alpha(M(h))&=&\alpha(h)\alpha(M(h')),\label{alt1}\\
\,\,\,\,\,\,\,\,\, \alpha(h')[N_{\beta}(h),x_{\alpha}]&=&\alpha(h)[N_{\beta}(h'),x_{\alpha}],
\forall \, h,h' \in H,\ \text{and}\ \alpha, \beta \in \Phi.\label{alt2}
\end{eqnarray}
Now, let $h \in H$ and $\alpha \in \Phi$, then
$$
\aligned
0&=[T(h),[x_{\alpha},x_{-\alpha}]]+[T(x_{\alpha}),[x_{-\alpha},h]]+[T(x_{-\alpha}),[h,x_{\alpha}]]
\\
&=
\sum_{\beta \in \Phi}[N_{\beta}(h),[x_{\alpha},x_{-\alpha}]]+\alpha(h)[T(x_{\alpha}),x_{-\alpha}]
+\alpha(h)[T(x_{-\alpha}),x_{\alpha}]
\\
&=
-\sum_{\beta \in \Phi}\beta([x_{\alpha},x_{-\alpha}])N_{\beta}(h)
+\alpha(h)[f_{\alpha}(x_{\alpha}),x_{-\alpha}]
+\alpha(h)\sum_{\beta \in \Phi}[g^{\alpha}_{\beta}(x_{\alpha}),x_{-\alpha}]
\\
&\quad 
+\alpha(h)[f_{-\alpha}(x_{-\alpha}),x_{\alpha}]
+\alpha(h)\sum_{\beta \in \Phi}[g^{-\alpha}_{\beta}(x_{-\alpha}),x_{\alpha}]
\\
& =
-\sum_{\beta \in \Phi}\beta([x_{\alpha},x_{-\alpha}])N_{\beta}(h)
-\alpha(h)\alpha(f_{\alpha}(x_{\alpha}))x_{-\alpha}
+\alpha(h)\sum_{\beta \in \Phi}[g^{\alpha}_{\beta}(x_{\alpha}),x_{-\alpha}]
\\
&\quad 
+\alpha(h)\alpha(f_{-\alpha}(x_{-\alpha}))x_{\alpha}
+\alpha(h)\sum_{\beta \in \Phi}[g^{-\alpha}_{\beta}(x_{-\alpha}),x_{\alpha}].
\endaligned
$$
The terms involving summations over $\beta\in\Phi$ can be rewritten as follows:
$$
\aligned
-\displaystyle{\sum_{\beta \in \Phi}}\beta([x_{\alpha},x_{-\alpha}])N_{\beta}(h)
&=
-\alpha([x_{\alpha},x_{-\alpha}])N_{\alpha}(h)+\alpha([x_{\alpha},x_{-\alpha}])N_{-\alpha}(h)\\
&\quad  -\displaystyle{\sum_{\beta \in \Phi - \{\pm \alpha\}}}\beta([x_{\alpha},x_{-\alpha}])N_{\beta}(h),
\\
\alpha(h)\sum_{\beta \in \Phi}[g^{\alpha}_{\beta}(x_{\alpha}),x_{-\alpha}]
& =
\alpha(h)[g^{\,\alpha}_{\,\alpha}(x_{\alpha}),x_{-\alpha}]+\alpha(h)\!\!\!\sum_{\beta \in \Phi - \{\alpha \}}\!\!\![g^{\alpha}_{\beta}(x_{\alpha})x_{-\alpha}],
\\
\alpha(h)\sum_{\beta \in \Phi}[g^{-\alpha}_{\beta}(x_{-\alpha}),x_{\alpha}]
& =
\alpha(h)[g^{-\alpha}_{-\alpha}(x_{-\alpha}),x_{\alpha}]
+\alpha(h)\!\!\!\sum_{\beta \in \Phi - \{-\alpha \}}\!\!\![g^{-\alpha}_{\beta}(x_{-\alpha})x_{\alpha}].
\endaligned
$$
Therefore,
$$
\aligned
0 &=
[T(h),[x_{\alpha},x_{-\alpha}]]+[T(x_{\alpha}),[x_{-\alpha},h]]+[T(x_{-\alpha}),[h,x_{\alpha}]]
\\
& =
-\alpha([x_{\alpha},x_{-\alpha}])N_{\alpha}(h)
+\alpha([x_{\alpha},x_{-\alpha}])N_{-\alpha}(h)-\!\!\!\!\!\sum_{\beta \in \Phi -{\{\pm \alpha}\}}\!\!\!\!\!\beta([x_{\alpha},x_{-\alpha}])N_{\beta}(h)
\\
&\quad
-\alpha(h)\alpha(f_{\alpha}(x_{\alpha}))x_{-\alpha}+\alpha(h)[g^{\alpha}_{\alpha}(x_{\alpha}),x_{-\alpha}]
+\alpha(h)\!\!\!\!\sum_{\beta \in \Phi - \{\alpha\}}\!\!\![g^{\alpha}_{\beta}(x_{\alpha}),x_{-\alpha}]
\\
&\quad
+\alpha(h)\alpha(f_{-\alpha}(x_{-\alpha}))x_{\alpha}\!+\!\alpha(h)[g^{-\alpha}_{-\alpha}(x_{-\alpha}),x_{\alpha}]
\!+\!\alpha(h)\!\!\!\!\!\sum_{\beta \in \Phi - \{-\alpha \}}\!\!\!\![g^{-\alpha}_{\beta}(x_{-\alpha}),x_{\alpha}].
\endaligned
$$
From this, it follows that:
$$
\aligned
-\alpha([x_{\alpha},x_{-\alpha}])N_{\alpha}(h)+\alpha(h)\alpha(f_{-\alpha}(x_{-\alpha}))x_{\alpha}
&=0,
\\
\alpha([x_{\alpha},x_{-\alpha}])N_{-\alpha}(h)-\alpha(h)\alpha(f_{\alpha}(x_{\alpha}))x_{-\alpha}
&=0,
\\
\alpha(h)[g^{\alpha}_{\alpha}(x_{\alpha}),x_{-\alpha}]+\alpha(h)[g^{-\alpha}_{-\alpha}(x_{-\alpha}),x_{\alpha}]
&=0,\quad \forall\, h \in H,\ \ \text{and}\ \ \alpha \in \Phi.
\endaligned
$$
Since $\alpha([x_{\alpha},x_{-\alpha}]) \neq 0$, we get,
\begin{equation}\label{alt3}
N_{\alpha}(h)=\frac{\alpha(f_{-\alpha}(x_{-\alpha}))}{\alpha([x_{\alpha},x_{-\alpha}])}\alpha(h)x_{\alpha},\,\forall h \in H,\, \forall \alpha \in \Phi.
\end{equation}
Since 
$\beta([x_{\beta},x_{-\beta}]) \neq 0$,
we may substitute \eqref{alt3} in \eqref{alt2}, to obtain,
$$
\alpha(h')\beta(h)\frac{\beta(f_{-\beta}(x_{-\beta}))}{\beta([x_{\beta},x_{-\beta}])}[x_{\beta},x_{\alpha}]=\alpha(h)\beta(h')\frac{\beta(f_{-\beta}(x_{-\beta}))}{\beta([x_{\beta},x_{-\beta}])}[x_{\beta},x_{\alpha}].
$$
This simplifies to,
\begin{equation}\label{alt4}
\alpha(h')\beta(h)\beta(f_{-\beta}(x_{-\beta}))[x_{\beta},x_{\alpha}]=\alpha(h)\beta(h')\beta(f_{-\beta}(x_{-\beta}))[x_{\beta},x_{\alpha}],
\end{equation}
for all $ h,h' \in H$ and for all $\alpha, \beta \in \Phi$. 
\smallskip
\noindent

We shall make a parenthesis here to prove the following standard result
where we use the hypothesis  that the rank of $\g$ is bigger than $1$.
\begin{Prop}\label{prop 1}{\sl
Let $\mathfrak{g}$ be a semi-simple Lie algebra of rank bigger than $1$. Fix 
a toral subalgebra $H\subset\g$ and let $\Phi$ be its corresponding root system.
Then, for each $\beta \in \Phi$, there exists a root
$\alpha \in \Phi$, such that $\beta+\alpha \in \Phi$ or $\beta-\alpha \in \Phi$.
}
\end{Prop}
\begin{proof} Let $(\,\cdot\,,\,\cdot\,)$ be the inner product defined in the 
$\mathbb{R}$-vector space generated by the roots in $\Phi$. 
If $\beta+\alpha \notin \Phi$ and $\beta-\alpha \notin \Phi$ for all $\alpha \in \Phi$, then, 
$(\beta,\alpha) \geq 0$ and $(\beta,\alpha) \leq 0$ for all $\alpha \in \Phi$
(see \cite{Hum}, {\bf Lemma 9.4}). Thus, $(\beta,\alpha)=0$ for all $\alpha \in \Phi$, 
which is a contradiction.
\end{proof}

\begin{Remark}\label{remark 6} 
We shall use this result as follows. If 
$\alpha,\beta \in \Phi$ are such that $\beta+\alpha \in \Phi$, 
by putting $\alpha^{\prime}=-\alpha \in \Phi$, we conclude that $\beta-\alpha^{\prime} \in \Phi$. 
So, we may always assume that for any given $\beta \in \Phi$, there are $\alpha,\alpha^{\prime} \in \Phi$, 
with $\alpha \neq \alpha^{\prime}$, such that $\beta+\alpha \in \Phi$ and $\beta-\alpha^{\prime} \in \Phi$.
\end{Remark}
\noindent
From now on, we shall assume that $\mathfrak{g}$ has rank bigger than $1$.

\smallskip
\noindent
Take $\alpha$ and $\beta$ in $\Phi$ so that $\alpha+\beta \in \Phi$. 
It follows from \eqref{alt4} that,
$$
\alpha(h')\beta(h)\beta(f_{-\beta}(x_{-\beta}))=\alpha(h)\beta(h')\beta(f_{-\beta}(x_{-\beta})),\quad\forall\, h,h' \in H.
$$
Suppose $\beta(f_{-\beta}(x_{-\beta})) \neq 0$. Then, $\alpha(h')\beta(h)=\alpha(h)\beta(h')$, for all $h,h' \in H$. In particular, for $h'=[x_{\alpha},x_{-\alpha}]$. Thus, $\beta=\frac{\beta([x_{\alpha},x_{-\alpha}])}{\alpha([x_{\alpha},x_{-\alpha}])}\alpha$. This means that $\beta$ is a root which is also a scalar multiple of the root $\alpha$,
thus implying that $\beta=\pm \alpha$. This however, contradicts the fact that $\alpha+\beta \in \Phi$. 
So, $\beta(f_{-\beta}(x_{-\beta}))=0$ for all $\beta \in \Phi$. It then follows from \eqref{alt3} that 
$N_{\alpha}=0$, for all $\alpha$, which in turn implies that $T(H) \subseteq H$.

\smallskip
\noindent
Now, consider $\alpha, \beta \in \Phi$, $\beta \neq  \pm \alpha$, and $h \in H$. Then,
\begin{multline}\label{alt5}
0=[T(h),[x_{\alpha},x_{\beta}]]+[T(x_{\alpha}),[x_{\beta},h]]+[T(x_{\beta}),[h,x_{\alpha}]]
\\
=(\alpha+\beta)(M(h))[x_{\alpha},x_{\beta}]-\beta(h)[T(x_{\alpha}),x_{\beta}]+\alpha(h)[T(x_{\beta}),x_{\alpha}]
\\
=(\alpha+\beta)(M(h))[x_{\alpha},x_{\beta}]-\beta(h)[f_{\alpha}(x_{\alpha}),x_{\beta}]
\\
- \beta(h)\sum_{\gamma \in \Phi}\,[g^{\alpha}_{\gamma}(x_{\alpha}),x_{\beta}]
+\alpha(h)[f_{\beta}(x_{\beta}),x_{\alpha}]
+\alpha(h)\sum_{\gamma \in \Phi}\,[g^{\beta}_{\gamma}(x_{\beta}),x_{\alpha}]
\\
=(\alpha+\beta)(M(h))[x_{\alpha},x_{\beta}]-\beta(h)\beta(f_{\alpha}(x_{\alpha}))x_{\beta}-\beta(h)[g^{\alpha}_{-\beta}(x_{\alpha}),x_{\beta}]
\\
-\beta(h)\!\!\!\!\sum_{\gamma \in \Phi -\{ -\beta\}}\!\!\!\!\![g^{\alpha}_{\gamma}(x_{\alpha}),x_{\beta}]+\alpha(h)\alpha(f_{\beta}(x_{\beta}))x_{\alpha}+\alpha(h)[g^{\beta}_{-\alpha}(x_{\beta}),x_{\alpha}]
\\
+\alpha(h)\!\!\!\!\sum_{\gamma \in \Phi -\{-\alpha \}}\!\!\!\!\![g^{\beta}_{\gamma}(x_{\beta}),x_{\alpha}].
\end{multline}
In this equation we may separate the terms belonging to $H$ in order to conclude that,
$$
\alpha(h)[g^{\beta}_{-\alpha}(x_{\beta}),x_{\alpha}]=\beta(h)[g^{\alpha}_{-\beta}(x_{\alpha}),x_{\beta}],\quad \forall h \in H.
$$
Since $g^{\beta}_{-\alpha}(x_{\beta}) \in \mathbb{F}x_{-\alpha}$ and $g^{\alpha}_{-\beta}(x_{\alpha}) \in \mathbb{F}x_{-\beta}$, there are scalars $\lambda^{\beta}_{-\alpha}, \lambda^{\alpha}_{-\beta} \in \F$, for which $g^{\beta}_{-\alpha}(x_{\beta})=\lambda^{\beta}_{-\alpha}x_{-\alpha}$ and $g^{\alpha}_{-\beta}(x_{\alpha})=\lambda^{\alpha}_{-\beta}x_{-\beta}$. Therefore, the last expression simplifies to,
$$
\alpha(h)\lambda^{\beta}_{-\alpha}[x_{-\alpha},x_{\alpha}]=\beta(h)\lambda^{\alpha}_{-\beta}[x_{-\beta},x_{\beta}],
\qquad \forall\,h \in H.
 $$
We now claim that $\lambda^{\beta}_{-\alpha}=\lambda^{\alpha}_{-\beta}=0$. We shall consider two cases: 
either $[x_{-\alpha},x_{\alpha}]$ and $[x_{-\beta},x_{\beta}]$ are linearly independent or they are not.

\smallskip
\noindent
\textbf{First case:} Assume $[x_{-\alpha},x_{\alpha}]$
and $[x_{-\beta},x_{\beta}]$ are linearly independent. The last expression implies that
$\alpha(h)\lambda^{\beta}_{-\alpha}=\beta(h)\lambda^{\alpha}_{-\beta}=0$, for all $h \in H$. 
Therefore, $\lambda^{\beta}_{-\alpha}=\lambda^{\alpha}_{-\beta}=0$.

\smallskip
\noindent
\textbf{Second case:} Assume $[x_{-\alpha},x_{\alpha}]$
and $[x_{-\beta},x_{\beta}]$ are linearly dependent.
Let $\nu \in \F$ be such that $[x_{-\beta},x_{\beta}]=\nu\,[x_{-\alpha},x_{\alpha}]$. Then, $\lambda^{\beta}_{-\alpha} \,\alpha=\nu \,\lambda^{\alpha}_{-\beta} \,\beta$. Since $[x_{-\alpha},x_{\alpha}]$ and $[x_{-\beta},x_{\beta}]$
are both different form $0$, we conclude that $\nu \neq 0$. This implies that $\beta$ is a non-zero saclar multiple of $\alpha$. Therefore, $\beta=\pm \alpha$, which is a contradiction.
It follows that if $\beta\ne \pm\alpha$, then 
$g^{\beta}_{-\alpha}(x_{\beta})=0$. Using this fact in \eqref{alt5}, we obtain,
$$
\aligned
0& =(\alpha+\beta)(M(h))[x_{\alpha},x_{\beta}]-\beta(h)\beta(f_{\alpha}(x_{\alpha}))x_{\beta}-\beta(h)[g^{\alpha}_{\alpha}(x_{\alpha}),x_{\beta}]\\
&\quad
-\beta(h)[g^{\alpha}_{-\alpha}(x_{\alpha}),x_{\beta}]+\alpha(h)\alpha(f_{\beta}(x_{\beta}))x_{\alpha}
+\alpha(h)[g^{\beta}_{\beta}(x_{\beta}),x_{\alpha}]\\
&\quad
+\alpha(h)[g^{\beta}_{-\beta}(x_{\beta}),x_{\alpha}].
\endaligned
$$
We may now separate the different terms of this expression that
belong to the different subspaces given by the root decomposition to conclude that,
\begin{eqnarray}
(\alpha+\beta)(M(h))[x_{\alpha},x_{\beta}]&=&
\beta(h)[g^{\alpha}_{\alpha}(x_{\alpha}),x_{\beta}]-\alpha(h)[g^{\beta}_{\beta}(x_{\beta}),x_{\alpha}],
\label{alt6}\\
\beta(h)[g^{\alpha}_{-\alpha}(x_{\alpha}),x_{\beta}]&=&\alpha(h)[g^{\beta}_{-\beta}(x_{\beta}),x_{\alpha}]=0,
\label{alt7}\\
\beta(h)\beta(f_{\alpha}(x_{\alpha}))&=&\alpha(h)\alpha(f_{\beta}(x_{\beta}))=0,\,\forall\,h \in H.\label{alt8}
\end{eqnarray}
Now, \eqref{alt7} implies that, 
$$
\lambda^{\alpha}_{-\alpha}\beta(h)[x_{-\alpha},x_{\beta}]
=\lambda^{\beta}_{-\beta}\alpha(h)[x_{-\beta},x_{\alpha}]=0,\quad \alpha,\beta \in \Phi,\ \ \text{and}\ \  h \in H.
$$
Since $\alpha\ne \pm\beta$, we may choose $\alpha$ and $\beta$ in $\Phi$ 
in such a way that $\alpha-\beta \in \Phi$ (see the {\bf Remark \ref{remark 6}} following {\bf Prop \ref{prop 1}}). 
Then, $\lambda^{\alpha}_{-\alpha}\beta=\lambda^{\beta}_{-\beta}\alpha =0$. 
Therefore, $\lambda^{\alpha}_{-\alpha} =\lambda^{\beta}_{-\beta} = 0$. 
Consequently, $g^{\alpha}_{-\alpha}=0$ for all $\alpha \in \Phi$.

\smallskip
\noindent
On the other hand, \eqref{alt8} implies that $\beta(f_{\alpha}(x_{\alpha}))=\alpha(f_{\beta}(x_{\beta}))=0$ for all $\beta \neq \pm \alpha$. But we had previously showed that
$\beta(f_{-\beta}(x_{-\beta}))=0$ for all $\beta \in \Phi$. 
And substituting $\beta$ by $-\beta$, we also get $\beta(f_{\beta}(x_{\beta}))=0$.
Therefore, $\beta(f_{\alpha}(x_{\alpha}))=0$ for all $\alpha,\beta \in \Phi$, 
which implies $f_{\alpha}=0$ for all $\alpha \in \Phi$. Thus,
so far we know that,
$$
T(x_{\alpha})=g^{\alpha}_{\alpha}(x_{\alpha})=\lambda^{\alpha}_{\alpha}x_{\alpha},\quad \alpha \in \Phi
\ \ \text{and}\ \ \lambda^{\alpha}_{\alpha} \in \mathbb{F}.
$$
Now, let $h \in H$ and $\alpha \in \Phi$. Then,
$$
\aligned
0
&=
[T(h),[x_{\alpha},x_{-\alpha}]]+[T(x_{\alpha}),[x_{-\alpha},h]]+[T(x_{-\alpha}),[h,x_{\alpha}]]
\\
&=\lambda^{\alpha}_{\alpha} \alpha(h)[x_{\alpha},x_{-\alpha}]
+\lambda^{-\alpha}_{-\alpha} \alpha(h)[x_{-\alpha},x_{\alpha}],\\
&=
(\lambda^{\alpha}_{\alpha}-\lambda^{-\alpha}_{-\alpha})\,
\alpha(h)\,[x_{\alpha},x_{-\alpha}],\quad h \in H\ \ \text{and}\ \ \alpha \in \Phi.
\endaligned
$$
It follows that $\lambda^{\alpha}_{\alpha}=\lambda^{-\alpha}_{-\alpha}$ for all $\alpha \in \Phi$. Thus, we write $\lambda_{\alpha}=\lambda^{\alpha}_{\alpha}=\lambda^{-\alpha}_{-\alpha}=\lambda_{-\alpha}$, and $T(x_{\alpha})=\lambda_{\alpha} \, x_{\alpha}$, for all $\alpha \in \Phi$.

\smallskip
\noindent 
Now let $\alpha,\beta,\gamma \in \Phi$. Then,
$$
\aligned
0&=
[T(x_{\alpha}),[x_{\beta},x_{\gamma}]]+[T(x_{\beta}),[x_{\gamma},x_{\alpha}]]
+[T(x_{\gamma}),[x_{\alpha},x_{\beta}]]\\
&=
\lambda_{\alpha}[x_{\alpha},[x_{\beta},x_{\gamma}]]
+\lambda_{\beta}[x_{\beta},[x_{\gamma},x_{\alpha}]]
+\lambda_{\gamma}[x_{\gamma},[x_{\alpha},x_{\beta}]].
\endaligned
$$
Using the fact that $\operatorname{ad}(x_\alpha)$ acts as
a $(1,1,1)$-derivation on $[x_\beta,x_\gamma]$, etc., we get,
$$
(\lambda_{\alpha}-\lambda_{\gamma})\left[[x_{\alpha},x_{\beta}],x_{\gamma}\right]
=(\lambda_{\beta}-\lambda_{\alpha})[x_{\beta},\left[x_{\alpha},x_{\gamma}]\right]
$$
Taking $\gamma=-\alpha$  and using the fact that $\lambda_{\alpha}=\lambda_{-\alpha}$ for all $\alpha \in \Phi$, this expression implies that $(\lambda_{\beta}-\lambda_{\alpha})\beta([x_{\alpha},x_{-\alpha}])x_{\beta}=0$, 
for all $\alpha,\beta \in \Phi$. 
Following \cite{Hum}, write $h_{\alpha}=[x_{\alpha},x_{-\alpha}]$ for all $\alpha \in \Phi$. 
Let $K_{\g}:\g \times \g \to \mathbb{F}$ be the Cartan-Killing form on $\g$. 
According to the {\bf Corollary} in {\bf \S 8.2} of \cite{Hum}, for each $\phi \in H^{*}$ 
there is a unique element $t_{\phi} \in H$
satisfying, $\phi(h)=K_{\g}(t_{\phi},h)$, for all $h \in H$. 
Moreover, the relationship between $h_{\alpha}$ and $t_{\alpha}$ is 
$h_{\alpha}=2\,t_{\alpha}/K_{\g}(t_{\alpha},t_{\alpha})$. 
On the other hand, the inner product $(\,\cdot\,,\,\cdot\,)$ defined in the 
$\mathbb{R}$-vector space generated by the roots is given by 
$(\beta,\alpha)=K_{\g}(t_{\beta},t_{\alpha})$, for all $\alpha,\beta$. 
Therefore, $\beta(h_{\alpha})=2(\beta,\alpha)/(\alpha,\alpha)$, for all $\alpha,\beta \in \Phi$. So, from $(\lambda_{\beta}-\lambda_{\alpha})\beta(h_{\alpha})=0$, we get 
$(\lambda_{\beta}-\lambda_{\alpha})(\beta,\alpha)=0$. This implies that
$\lambda_{\alpha}=\lambda_{\beta}$, whenever $(\beta,\alpha) \neq 0$.
Now, from {\bf Lemmas A} and {\bf D} in {\bf \S10.4} of \cite{Hum} one knows that there is
a unique maximal root $\beta^{\prime}$, such that, 
$(\beta^\prime,\beta^\prime) \geq (\beta^\prime,\alpha) \geq (\alpha,\alpha)$ for all $\alpha \in \Phi$. 
If $\alpha$ satisfies $(\beta^{\prime},\alpha)=0$, then $(\alpha,\alpha)=0$, 
which contradicts the fact that $(\,\cdot\,,\,\cdot\,)$ is positive definite.
Therefore, $(\beta^{\prime},\alpha) \neq 0$, for all $\alpha \in \Phi$.
It follows that $\lambda_{\beta^{\prime}}=\lambda_{\alpha}$ for all $\alpha \in \Phi$. 
Thus, we may write $\lambda=\lambda_{\alpha}=\lambda_{\beta^{\prime}}$;
and since the scalar does not depend on the roots we have,
$T(x_{\alpha})=\lambda\,x_{\alpha}$ for all $\alpha \in \Phi$.

\smallskip
\noindent 
Let $\alpha \in \Phi$. It follows from \eqref{alt1} that, $\alpha \circ M=\frac{\alpha(M([x_{\alpha},x_{-\alpha}]))}{\alpha([x_{\alpha},x_{-\alpha}])}\,\alpha$.
We write $c_{\alpha}=\frac{\alpha(M([x_{\alpha},x_{-\alpha}]))}{\alpha([x_{\alpha},x_{-\alpha}])}\,\alpha$. Let $\alpha,\beta \in \Phi$ be such that $\alpha+\beta \in \Phi$. From \eqref{alt6} we get, 
$(\alpha+\beta)M=\lambda \alpha+\lambda \beta$. 
Since, $\gamma \circ M=c_{\gamma}\gamma$ for all $\gamma \in \Phi$, 
we get $(\lambda-c_{\alpha})\alpha=(\lambda-c_{\beta})\beta$. Thus, 
$\alpha+\beta\in\Phi$ implies that
$c_{\alpha}=\lambda=c_{\beta}$. 
Since $\alpha$ is arbitrary, we conclude that $\alpha \circ M=\lambda \,\alpha$, for any $\alpha \in \Phi$.
That is, $M^\ast(\alpha)=\lambda\,\alpha$ for any root. Since the roots span $H^\ast$
and $M=\pi_H\,\circ\,(T\vert_H):H\to H$, it follows that $M=\lambda\,\alpha$. 
Consequently, $T=\lambda \operatorname{Id}_{\mathfrak{g}}$.
\medskip

{\boldmath $(2)$} Let $T \in \operatorname{HL}(\mathfrak{sl}_2(\F))$. 
Then, the matrix of $T$ in the basis $\{H,E,F\}$, has the form,
$$
T=\begin{pmatrix}
T_{11} & T_{12} & T_{13} \\
2 T_{13} & T_{22} & T_{23} \\
2 T_{12} & T_{32} & T_{33}
\end{pmatrix}\!\!,\quad T_{i,j} \in \mathbb{F}.
$$
Then,
$$
T=
\begin{pmatrix}
\displaystyle{\frac{2(T_{11}-T_{22})}{3}}
& T_{12} & T_{13} \\
2 T_{13} & \displaystyle{\frac{-T_{11}+T_{22}}{3}} & T_{23} \\
2 T_{12} & T_{32} & \displaystyle{\frac{-T_{11}+T_{22}}{3}}
\end{pmatrix}
+\frac{1}{3}(T_{11}+2T_{22})\operatorname{Id}_{3 \times 3}.
$$
This proves that $\operatorname{HL}(\mathfrak{sl}_2(\F))=\operatorname{SHL}(\mathfrak{sl}_2(\F)) \oplus \F \operatorname{Id}_{\mathfrak{sl}_2(\F)}$, where 
$$
\operatorname{SHL}(\mathfrak{sl}_2(\F))=\{T \in \operatorname{HL}(\mathfrak{sl}_2(\F))\,|\,\operatorname{Tr}(T)=0\}.
$$
Observe that $\operatorname{SHL}(\mathfrak{sl}_2(\F))$ is equal to $\Der_{(-1,1,1)}(\mathfrak{sl}_2(\F))$. The $\mathfrak{sl}_2(\F)$-action on $\operatorname{HL}(\mathfrak{sl}_2(\F))$ is given by,
$$
\aligned
x.T&=[\ad_{\mathfrak{sl}_2(\F)}(x),T]_{\gl(\mathfrak{sl}_2(\F))}=\ad_{\mathfrak{sl}_2(\F)}(x) \circ T-T \circ \ad_{\mathfrak{sl}_2(\F)}(x), 
\\
(x.T)(y)&=[x,T(y)]-T([x,y]),\quad \forall \, T \in \operatorname{HL}(\mathfrak{sl}_2(\F)),\quad \forall \, x,y \in \mathfrak{sl}_2(\F).
\endaligned
$$
Since $\Tr(x.T)=0$, for all $x \in \mathfrak{sl}_2(\F)$ and $T \in \operatorname{HL}(\mathfrak{sl}_2(\F))$, 
it follows that, $\Der_{(-1,1,1)}(\mathfrak{sl}_2(\F))$ is an $\mathfrak{sl}_2(\F)$-submodule of $\operatorname{HL}(\mathfrak{sl}_2(\F))$ such that $x.\operatorname{HL}(\mathfrak{sl}_2(\F)) \subset \Der_{(-1,1,1)}(\mathfrak{sl}_2(\F))$, for all $x \in \mathfrak{sl}_2(\F)$. On the other hand,
$$
\aligned
H.\begin{pmatrix}
T_{11} & T_{12} & T_{13} \\
2 T_{12} & T_{22} & T_{23} \\
2 T_{13} & T_{32} & T_{33}
\end{pmatrix}&=
\begin{pmatrix}
0 & -2 T_{12} & 2 T_{13} \\
4 T_{13} & 0 & 4 T_{23} \\
-4 T_{12} & -4 T_{32} & 0
\end{pmatrix}\!,\\
E.\begin{pmatrix}
T_{11} & T_{12} & T_{13} \\
2 T_{13} & T_{22} & T_{23} \\
2 T_{12} & T_{32} & T_{33}
\end{pmatrix}&=
\begin{pmatrix}
4 T_{12} & T_{32} & -T_{11}+T_{22} \\
2(-T_{11}+T_{22}) & -2T_{12} & -T_{13} \\
2 T_{32} & 0 & -2 T_{12}
\end{pmatrix}\!,\\
F.\begin{pmatrix}
T_{11} & T_{12} & T_{13} \\
2 T_{13} & T_{22} & T_{23} \\
2 T_{12} & T_{32} & T_{33}
\end{pmatrix}&=
\begin{pmatrix}
-4T_{13} & T_{11}-T_{22} & -T_{23} \\
-2 T_{23} & 2 T_{13} & 0 \\
2(T_{11}-T_{22}) & 4 T_{12} & 2 T_{13}
\end{pmatrix}\!.
\endaligned
$$
Let $\xi \in \mathbb{F}$, and write $\operatorname{HL}(\mathfrak{sl}_2(\F))_{\xi}:=\{T \in V\,|\,H.T=\xi T\}$. 
The right-hand sides are clearly written in terms of the basis
$\{P,Q,R,S,T\}$ given by, $P=\left(\begin{smallmatrix}
2 & 0 & 0\\
0 & -1 & 0\\
0 & 0 & -1
\end{smallmatrix}\right)$, $Q=\left(\begin{smallmatrix}
0 & 1 & 0\\
0 & 0 & 0\\
2 & 0 & 0
\end{smallmatrix}\right)$, $R=\left(\begin{smallmatrix}
0 & 0 & 1\\
2 & 0 & 0\\
0 & 0 & 0
\end{smallmatrix}\right)$, $S=
\left(\begin{smallmatrix}
0 & 0 & 0\\
0 & 0 & 1\\
0 & 0 & 0
\end{smallmatrix}\right)$ and $T=
\left(\begin{smallmatrix}
0 & 0 & 0\\
0 & 0 & 0\\
0 & 1 & 0
\end{smallmatrix}\right)$.
In fact,
$$
\operatorname{HL}(\mathfrak{sl}_2(\F))_{4}=\F \,S,\quad 
\operatorname{HL}(\mathfrak{sl}_2(\F))_2=\F \,R,\quad 
\operatorname{HL}(\mathfrak{sl}_2(\F))_0=\F \,P\,\oplus \F \operatorname{Id}_{3 \times 3},
$$
$$
\operatorname{HL}(\mathfrak{sl}_2(\F))_{-2}=\F \,Q,\quad
\operatorname{HL}(\mathfrak{sl}_2(\F))_{-4}=\F\, T.
$$
From the above, it is clear that $\operatorname{SHL}(\mathfrak{sl}_2(\F))$ is an irreducible $\mathfrak{sl}_2(\F)$-submodule, because $\operatorname{SHL}(\mathfrak{sl}_2(\F))_0$ is one dimensional.
\end{proof}

\noindent
\begin{Remark}{\rm
It is only natural to ask what difference would it make if the
hypotheses are relaxed so as to let $\g$ be a semisimple Lie algebra.
In the course of the last proof above we were faced at some point with the fact that
$(\lambda_{\beta}-\lambda_{\alpha})(\beta,\alpha)=0$, implied that
$\lambda_\beta=\lambda_\alpha$, whenever $(\beta,\alpha)\ne 0$.
This conclusion is no longer true if $\g$ is semisimple  because there
may be several independent simple components on which
$\alpha$ and $\beta$ lie, and still have $(\beta,\alpha)=0$.}
\end{Remark}

\section*{Acknowledgements}
The author thanks the support 
provided by the post-doctoral fellowships
FOMIX-YUC 221183, FORDECYT 265667 and CONACYT 000153
that allowed him to complete this work
at CIMAT - Unidad M\'erida. 
The author would also like to deeply thank Professors G. Salgado and O. A. S\'anchez-Valenzuela
for suggesting some of the problems addressed in this work, reviewing the manuscript and 
making valuable suggestions that greatly improved the final version of this paper. 

\end{document}